\documentclass[a4paper,12pt]{amsart}

\usepackage[utf8]{inputenc}
\usepackage[T1]{fontenc}
\usepackage{amssymb}
\usepackage[pagebackref]{hyperref}
\usepackage{bbm}
\usepackage{mathtools}
\usepackage{enumerate}
\usepackage{ulem}
\usepackage[left=2cm,right=2cm]{geometry}
\allowdisplaybreaks
\newtheorem{theorem}{Theorem}[section]
\newtheorem{lemma}[theorem]{Lemma}
\newtheorem{cor}[theorem]{Corollary} 
\newtheorem{prop}[theorem]{Proposition} 

\theoremstyle{definition}

\theoremstyle{remark}
\newtheorem{remark}[theorem]{Remark}

\numberwithin{equation}{section}

\def\ii{{\bf i}}

\newcommand{\id}{\mathrm{Id}}
\newcommand{\ind}[1]{\mathbf{1}_{\{#1\}}}

\newcommand{\E}{\mathbb{E}} 
\newcommand{\p}{\mathbb{P}} 
\newcommand{\R}{\mathbb{R}} 
\newcommand{\C}{\mathbb{C}}

\newcommand{\N}{\mathbb{N}}
\newcommand{\XX}{\mathbf{X}}

\newcommand{\xX}{\mathbb{X}}
\newcommand{\xx}{\mathbf{x}}
\newcommand{\yy}{\mathbf{y}}
\newcommand{\yY}{\mathbb{Y}}
\newcommand{\hxx}{\widehat{\mathbf{x}}}
\newcommand{\hw}{\widehat{w}}
\newcommand{\hx}{\widehat{x}}
\newcommand{\hb}{\widehat{b}}
\newcommand{\hz}{\widehat{z}}

\newcommand{\YY}{\mathbf{Y}}

\newcommand{\1}{\mathbf{1}}
\DeclareMathOperator{\Var}{Var}	
\DeclareMathOperator{\Cov}{Cov}
\DeclareMathOperator{\Med}{Med}	

\DeclareMathOperator{\dist}{dist}
\DeclareMathOperator{\Ent}{Ent}
\DeclareMathOperator{\tr}{tr}

\title[]{Concentration of measure for non-linear random matrices with applications to neural networks and non-commutative polynomials}

\author{Radosław Adamczak}
\address{University of Warsaw, Institute of Mathematics, ul. Banacha 2, 02-097 Warszawa, Poland}
\email{r.adamczak@mimuw.edu.pl}
\thanks{Research partially supported by The National Science Center, Poland, grant IMPRESS-U 2023/05/Y/ST1/00188}

\date{\today}

\begin{document}
\maketitle
\begin{abstract}We prove concentration inequalities for several models of non-linear random matrices. As corollaries we obtain estimates for linear spectral statistics of the conjugate kernel of neural networks and non-commutative polynomials in (possibly dependent) random matrices.
\end{abstract}

\section{Introduction}

The theory of global behaviour of eigenvalues of classical ensembles of random matrices, such as Wigner matrices or sample covariance matrices, is by now well understood. In particular, the foundational results due to Wigner \cite{MR95527} and Marchenko--Pastur \cite{MR208649}, concerning convergence of the empirical spectral distribution have been complemented with central limit theorems (e.g., \cite{MR2561434,MR2556016}), large deviations (e.g., \cite{MR1465640,MR1743092,MR3265172,MR4869057}), and concentration inequalities (e.g., \cite{MR1781846,MR2737712,MR2535081,MR3055255,MR4049090}), the list of references being clearly non-exhaustive.

In recent years, with rapidly increasing interest in neural networks, new models of random matrices have attracted considerable attention, with hope that investigation of their spectral properties may give important information concerning the behaviour of learning algorithms at initialization or at an early learning stage. In particular, Pennington and Worah \cite{MR4063603}, Benigni and P\'{e}ch\'{e} \cite{MR4346666}, and Guionnet and Piccolo \cite{guionnet2024global} obtained limit theorems for the spectral measure of the conjugate kernel of neural networks at initialization as the size of the data sample and the width of the network layers go to infinity at the same speed. Due to the presence of non-linear activation function, applied to the entries of a matrix which evolves along the network, the corresponding models are sometimes called \emph{non-linear random matrices}.

The goal of this paper is to investigate concentration properties of such matrices, in particular to obtain concentration inequalities for their spectral statistics. We obtain several concentration estimates for Lipschitz functions, with strength depending on the properties of the activation functions and the distribution of the network. Inequalities for spectral statistics are then obtained as a corollary to a more general form of concentration of measure, similarly as in the by now classical approach by Guionnet--Zeitouni \cite{MR1781846} in the case of Wigner matrices. From the perspective of neural networks, non-asymptotic concentration estimates are of interest, because they may be useful in the analysis of the learning process \cite{arous2023highdimensionalsgdalignsemerging} or to quantify the learning error in special classes of models \cite{MR4522365}. We investigate the concentration phenomenon at initialization, however it was recently proved in \cite{wang2024spectral} that the spectral distribution of a neural network remains stable during the initial stage of the gradient descent evolution. Moreover, certain properties of random matrices related to the neural network at initialization are crucial elements in the analysis of the speed of convergence of optimization algorithms used for learning, see, e.g., \cite{MR4412181} and the references therein.

Since most of our proofs allow for some amount of dependence, we do not restrict the results to the usual setting of neural networks in which consecutive layers are initialized independently. An additional motivation for a more general approach comes from the fact that inequalities we obtain in the dependent setting are often equally strong as those in the independent one, moreover the method of proof allows to obtain analogous results for non-commutative polynomials in random matrices, leading to non-asymptotic estimates for their general spectral statistics. This line of research is related to the theory of free probability, allowing in particular to reduce almost sure freeness  to freeness in expectation (we refer, e.g., to the monographs \cite{MR3585560,MR1746976} for an exposition of free probability). Prior results in this direction were obtained by Meckes and Szarek \cite{MR2869165} and more recently by Bao and Munoz George \cite{bao2024ultrahighordercumulants} who investigated the special case of the trace.

\medskip
The organization of the article is as follows. In Section \ref{sec:notation} we introduce the basic setting, in particular the assumptions on the random matrices we consider, and we illustrate them with some examples. We also state some simple facts related to the general theory of concentration of measure. In Section \ref{sec:results} we first state our results, dividing them according to the assumptions on the properties of the random matrix ensemble. Subsequently, we present various applications of concentration inequalities to the theory of random matrices (Section \ref{sec:applications}). Finally, we provide a discussion of optimality and necessity of various assumptions (Section \ref{sec:discussion}). Section \ref{sec:proofs} is devoted to the proofs of our main theorems and Section \ref{sec:applications-proofs} to the proofs of corollaries. In the Appendix we justify some simple technical lemmas which are used in our arguments.

\section{Notation and setting}\label{sec:notation}

\subsection{Preliminary notation} Unless specified otherwise, in this article, for a function $\sigma \colon \R \to \R$ and for an $m_1\times m_2$ matrix $\xx = (x_{ij})_{1\le i\le m_1,1\le j\le m_2}$, by $\sigma(\xx)$ we will denote the matrix $(\sigma(x_{ij}))_{1\le i\le m_1,1\le j\le m_2}$, i.e., we will apply the function $\sigma$ entry-wise. This is in contrast with the \emph{spectral convention} used usually in Random Matrix Theory but agrees with the notation used in the neural networks literature. To avoid confusion, when speaking about spectral statistics of a matrix we will define them explicitly in terms of eigenvalues or singular values.

For a vector $x = (x_1,\ldots,x_n) \in \R^n$  let $\|x\|_2 = (\sum_{i=1}^n x_i^2)^{1/2}$ be its standard Euclidean norm.
The Hilbert--Schmidt (or Frobenius) norm of an $m_1 \times m_2$ matrix $\xx = (x_{ij})_{1\le i\le m_1,1\le j\le m_2}$ is defined as $\|\xx\|_{HS} = \sqrt{\sum_{i=1}^{m_1}\sum_{j=1}^{m_2} x_{ij}^2}$. By $\|\xx\|_{op}$ we will denote the operator norm, i.e., $\|x\|_{op} = \sup_{y \in \R^{m_2}, \|y\|_2 \le 1} \|\xx y\|_2$.
We will often identify spaces of $m_1 \times m_2$ matrices with $\R^{m_1m_2}$, sometimes writing $\R^{m_1\times m_2}$ to stress the underlying matrix structure. As we will encounter also products of such spaces (since neural networks are defined in terms of several matrices), we will often identify, e.g. $\R^{m_1\times m_2}\times \R^{m_3\times m_4} \times \R^{m_5}$ with $\R^{m_1m_2+m_3m_4 + m_5}$ without making additional comments. Such identifications should always be clear from the context. In such situations, for $\xx \in \R^{m_1\times m_2}\times \R^{m_3\times m_4} \times \R^{m_5}$, by $\|\xx\|_2$ we mean the standard Euclidean norm of $\xx$ as an element of  $\R^{m_1m_2+m_3m_4 + m_5}$.

For a function $f \colon \R \to \R$, by $\|f\|_\infty$ we denote the sup norm, i.e. $\|f\|_\infty = \sup_{x \in \R} |f(x)|$.

We will usually denote random matrices by capital letters and deterministic matrices by small letters. Since quite often we will discuss sequences of matrices, we may occasionally write, e.g., $W_k(i,j)$ to denote the $(i,j)$ entry of the matrix $W_k$.

For a positive integer $n$, by $[n]$ we will denote the set $\{1,\ldots,n\}$.

The rest of our notation is rather standard and, if needed, it will be introduced in subsequent parts of the article.

\subsection{Main setting}
Let $L$, $n_0,\ldots,n_{L}$ be positive integers and consider $W_\ell,B_\ell$, $\ell=0,\ldots,L-1$, where $W_\ell$ is an $n_{\ell+1} \times n_{\ell}$ random matrix, and $B_\ell$ a random vector in $\R^{n_{\ell+1}}$ (we will treat all vectors as \emph{vertical} ones). Consider also a sequence of activation functions $\sigma_\ell \colon \R \to \R$, $\ell=0,\ldots,L-1$. Let $X$ be a random matrix of size $n_0 \times n$. Define recursively the matrices $Z_\ell$, $\ell = 0,\ldots,L$ as $Z_0 = \frac{1}{\sqrt{n_0}}X$,
\begin{align}\label{eq:model}
  Z_{\ell+1} = \frac{1}{\sqrt{n_{\ell+1}}}\sigma_\ell\Big(W_\ell Z_\ell + B_\ell \1_{n}^T \Big),
\end{align}
where $\1_{n}$ denotes the $n$-dimensional vector with all coordinates equal to 1. Thus $Z_\ell$ is an $n_{\ell} \times n$ random matrix.

In terms of neural networks, we may think of the columns of $X$ as the $n$ samples of $n_0$-dimensional data and $n_1,\ldots,n_L$ as the widths of consecutive layers of a neural network with weight matrices $W_0,\ldots,W_{L-1}$ and bias vectors $B_0,\ldots,B_{L-1}$. The $j$-th column of $Z_L$ is the result of the action of the network on the $j$-th column of $X$. We refer, e.g., to \cite{MR4719738,MR4604111,MR3617773} for a general introduction of the theory of neural networks.

We will be interested in concentration properties of random variables of the form $F(Z)$ where $F \colon \R^{n_L \times n} \to \R$ is Lipschitz with respect to the Hilbert--Schmidt norm. The primary example we have in mind are the statistics of singular values of $Z_L$, i.e.,
\begin{align}\label{eq:les}
  S_f(Z_L) = \sum_{i=1}^n f(\gamma_i),
\end{align}
where $\gamma_i$ are the singular values of $Z_L$, i.e., the square roots of the eigenvalues of $Z_L^T Z_L$, and $f\colon \R\to \R$. The Hoffman--Wielandt inequality \cite{MR52379} (see Theorem \ref{thm:HW} below) implies that if $f$ is $1$-Lipschitz, then $S_f \colon \R^{n_L \times n} \to \R$ is $\sqrt{n}$-Lipschitz. Another example of interest is the operator norm of $Z_L$ (i.e., the largest singular value), which is 1-Lipschitz.

Let us remark that $Z_L^TZ_L$ is the matrix known in the neural network setting as the conjugate kernel. It is quite important in applications as it describes the inner products between features extracted from the data by the network. See, e.g., \cite{MR4063603,MR4346666,guionnet2024global,pmlr-v139-hu21b} for results concerning the limiting distribution of the empirical spectral measure of the conjugate kernel (we refer to Section \ref{sec:applications} for the definition of the spectral measure) and, e.g., \cite{howard2024the,zeng2024featuremappingphysicsinformedneural} for some recent applications of the conjugate kernel in training neural networks. In what follows we will also address concentration properties for $S_f(Z_L^TZ_L) = \sum_{i=1}^n f(\gamma_i^2)$ and $\|Z_L^TZ_L\|_{op}$. We remark that in some of the references discussed above the authors consider the matrix $Z_LZ_L^T$ and not the conjugate kernel $Z_L^TZ_L$, however it is straightforward to obtain the spectral distribution of one of them from the other one.

\begin{remark}
Let us briefly describe the motivation behind our normalisation by $\sqrt{n_L}$. It ensures that in a typical situation the entries of $Z_L$ have the second moment of the order $1/n_L$, and as a consequence, $\frac{1}{n} \E \sum_{i=1}^n \gamma_i^2 = \frac{1}{n}\E \|Z_L\|_{HS}^2$ is of the order of a constant. In other words the second order of the expected spectral measure of $Z_L$ is of constant order which ensures that the sequence of expected spectral measures is tight when $n \to \infty$. This is a typical normalisation in Random Matrix Theory, allowing to study the global (bulk) behaviour of singular values.

At the same time $\frac{1}{n_\ell}$ is the order of variance used in the most common initializations of weights in random matrices (such as LeCun \cite{LeCun1998} or He \cite{7410480} initialization and also Xavier initialization \cite{pmlr-v9-glorot10a} if all layers have comparable width). Here we prefer to assume that the matrices $W_i$ have entries of constant order and hide this scaling factor in $Z_\ell$, which allows for a simpler formulation of our concentration assumptions below.
\end{remark}

\subsection*{Concentration of measure assumptions}
We will say that a random vector $\XX$ with values in $\R^N$ satisfies the subgaussian concentration property with constant $\rho$ if for any $F \colon \R^N \to \R$, 1-Lipschitz with respect to the standard Euclidean norm and any $t>0$,
\begin{align}\label{eq:concentration-assumption}
  \p(|F(\XX) - \Med F(\XX)| \ge t) \le 2\exp(-t^2/\rho^2)
\end{align}
where $\Med F(\XX)$ is (any) median of $F(\XX)$.

In our results we will assume that this condition holds for $\XX = (X,(W_\ell,B_\ell)_{\ell=0,\ldots,L-1})$ seen as a random vector in $\R^{N}$, with $N = nn_0 + \sum_{\ell=0}^{L-1} (n_\ell n_{\ell+1} + n_{\ell+1})$.

In most of the article we will consider real matrices. One exception concerns Section \ref{sec:polynomials}. In this case we will identify $\C^M$ with $\R^{N}$ for $N = 2M$ without making further comments.

\subsection*{Examples}

Let us provide a few examples of random vectors satisfying the subgaussian concentration property.

\subsubsection*{The log-Sobolev inequality} The inequality \eqref{eq:concentration-assumption} holds if $\XX$ satisfies the log-Sobolev inequality
\begin{align}\label{eq:LSI}
  \Ent f^2(X) \le C \E |\nabla f(X)|^2
\end{align}
for all sufficiently smooth functions. Here, for a non-negative random variable $Y$, by $\Ent Y$ we denote the entropy of $Y$, given by $\Ent Y = \E Y\ln Y - \E Y\ln \E Y$.

In this case $\rho$ depends only on $C$. It is well known that the log-Sobolev inequality tensorizes, i.e., if $X_1,\ldots,X_n$ are independent and satisfy \eqref{eq:LSI} then so does $\XX = (X_1,\ldots,X_n)$. In particular, if the components of $\XX$ are independent and satisfy \eqref{eq:LSI}, then $\XX$ satisfies the subgaussian concentration property.

This observation was used by Guionnet and Zeitouni to obtain concentration inequalities for linear eigenvalue statistics of Wigner matrices. In our case it implies that the data $\XX = (X,(W_\ell,B_\ell)_{\ell=0,\ldots,L-1})$ satisfies the subgaussian concentration property if all the entries of $X$, $W_\ell$, $B_\ell$ are independent and satisfy the one-dimensional log-Sobolev inequality. We remark that in the case of one-dimensional random variables, Bobkov and G\"otze \cite{MR1682772} found effective necessary and sufficient conditions for \eqref{eq:LSI} to hold. It is well known that this inequality is satisfied, e.g., by standard Gaussian random variables or random variables uniform on $(-1,1)$, which are often used for initialization of neural networks.

Beyond the product case, examples of measures satisfying the LSI are provided, e.g., by the famous Bakry--\'{E}mery criterion \cite{MR889476}. The log-Sobolev inequality holds if $\XX$ has a density with respect to the Lebesgue measure of the form $e^{-V}$ where $D^2V \ge \lambda \id$ for some $\lambda > 0$. It is also known that the log-Sobolev inequality is preserved (with changed constants) if one replaces the law $\mu_\XX$ of $\XX$ by a measure possessing a density with respect to $\mu_\XX$ which is separated from zero and infinity (see, e.g., \cite[Th\'eor\`eme 3.4.3]{MR1845806}).

Other examples related to this context, which are of interest from the random matrix point of view, are multi-matrix models, studied, e.g., by Figalli and Guionnet \cite{MR3646880} or Guionnet and Maurel-Segala \cite{MR2353386}. These are models of $L$ matrices $X_1,\ldots,X_L$ of size $n\times n$, with joint density of the form
\begin{align}\label{eq:multi}
  \p((X_1,\ldots,X_L) \in dx_1\cdots dx_L)
  = \frac{1}{\mathcal{Z}_n} e^{n \tr V(n^{-1/2}x_1,\ldots,n^{-1/2} x_L)}\prod_{\ell=1}^n \ind{\|x_\ell\|_{op} \le M\sqrt{n}}dx_1\cdots dx_L,
\end{align}
where $V$ is a non-commutative polynomial (we allow here the case $M = \infty$) and $\mathcal{Z}_n$ is the normalizing constant. We do not want to enter into a detailed technical discussion, but we point out that the articles \cite{MR3646880, MR2353386} present conditions (based on the Bakry-\'{E}mery theory) on the polynomial $V$ (and possibly the truncation level $M$) which guarantee that such ensembles satisfy the log-Sobolev inequality with a constant independent of the dimension $n$ (we remark that our normalization of $X_\ell$'s differs by a factor $\sqrt{n}$ with respect to the one used in \cite{MR3646880, MR2353386}, we chose this normalization to make the example consistent with our general notation). In particular our results will apply to non-commutative polynomials in the random matrices $X_1,\ldots,X_n$.

\subsubsection*{Uniform measure on the sphere and Haar-distributed orthogonal matrices} The subgaussian concentration property with dimension independent constants holds if $\XX$ is a random vector uniformly distributed on $\sqrt{N}S^{N-1}$, where $S^{N-1}$ is the unit sphere in $\R^N$. It also holds if $\XX = \sqrt{N} U$, where $U$ is a Haar-distributed random orthogonal matrix. Such initialization distributions for $B_\ell$ and $W_\ell$ with dependencies among entries have also been considered in the literature on neural networks \cite{10.5555/3295222.3295232,pennington2018emergence,pastur2022eigenvalue}. In fact these two examples are also related to the log-Sobolev inequality and the Bakry--\'{E}mery condition but on the respective manifolds. Concentration follows from comparison of the geodesic and Euclidean distances (see, e.g., \cite{MR2760897,MR3837270}).

\subsubsection*{Tensorization of concentration} A condition guaranteeing subgaussian concentration, weaker than the log-Sobolev inequality is Talagrand's transport-entropy inequality
\begin{displaymath}
\mathcal{W}_2^2(\nu,\mu) \le C H(\nu|\mu) \textrm{ for all probability measures $\nu$ on $\R^N$},
\end{displaymath}
where $\mu$ is the law of a random vector $\XX$. Recall that for $p\ge 1$,  $\mathcal{W}_p$ is the $p$-th  Wasserstein distance, defined by
\begin{equation}\label{eq:Wasserstein}
\mathcal{W}_p^p(\nu,\mu) = \inf_{(Y,Z)}\E\|Y-Z\|_2^p,
\end{equation}
where the  infimum is taken over all couplings $(Y,Z)$ of $\nu,\mu$, i.e., pairs of random variables $(Y,Z)$ defined on a common probability space, such that $Y$ has the law $\nu$ and $Z$ has the law $\mu$, while $H(\nu|\mu)$ is the Kullback-Leibler divergence defined as
\begin{displaymath}
  H(\nu|\mu) = \int_{\R^N} \ln \Big(\frac{d\nu}{d\mu}\Big)d\nu
\end{displaymath}
if $\nu$ is absolutely continuous with respect to $\mu$, and $\infty$ otherwise.

It has been proved by Gozlan \cite{MR2573565} that if $X,X_1,X_2,\ldots$ are i.i.d. random vectors, then $\XX_L = (X_1,\ldots,X_L)$ satisfies the subgaussian concentration inequality for all $L=1,2,\ldots$ with a constant independent of $L$ if and only if $X$ satisfies Talagrand's inequality, while an earlier result by Bobkov and G\"otze asserts that in fixed dimension subgaussian concentration is equivalent to a related transport-entropy inequality involving the first Wasserstein distance (we refer to \cite{MR1682772} for details). Similarly as with the log-Sobolev inequality, there are characterizations of transport-entropy inequalities for real random variables \cite{MR2352486}. They imply in particular that subgaussian concentration is a much weaker property than dimension-free subgaussian concentration. However, since in our applications we will allow the constants to depend on $L$ we may take advantage of the following easy observation.
\begin{prop}\label{prop:independent-concentration}
Let $X_i$ be a random vector with values in $\R^{n_i}$, $i = 1,\ldots,L$. Assume that $X_i$'s are independent. If all the vectors $X_i$ satisfy the subgaussian concentration property with constant $\rho$, then $\XX = (X_1,\ldots,X_L)$ satisfies the subgaussian concentration property with constant $C\sqrt{L}\rho$, where $C$ is a universal constant.
\end{prop}

The above proposition is a part of folklore. For completeness we will provide its proof in the appendix.

\subsubsection*{Preservation of concentration properties}

Let us finally observe that if $\XX$ satisfies the subgaussian concentration property, then clearly so does $\YY = T(X)$, where $T$ is any 1-Lipschitz map. An easy consequence of this observation is the following proposition.

\begin{prop}\label{prop:concentration-repetition}If $(X_1,\ldots,X_L)$ satisfies the subgaussian concentration property with constant $\rho$ and $\YY = (X_{i_1},\ldots,X_{i_m})$ for some $i_1,\ldots,i_m \in [L]$, then $\YY$ satisfies the subgaussian concentration property with constant $\rho \sqrt{R}$, where $R = \max_{\ell \in [L]} |\{j\in [m]\colon i_j = \ell\}|$.
\end{prop}
\begin{remark}
It is well known that subgaussian concentration around the mean and median are equivalent, more precisely, if one of them holds with constant $\rho$, then the other one holds with constant $C\rho$, where $C$ is a  universal constant (see Lemmas \ref{le:concentration-around-mean} and \ref{le:concentration-equivalence} below). We choose to state our assumptions in terms of concentration around the median, since this will be more convenient in some of the proofs, however we will often jump back and forth between the two formulations.
\end{remark}

\section{Main results}\label{sec:results}

In the setting of nonlinear random matrices $Z_L$, as defined in \eqref{eq:model}, we will be interested mostly in concentration properties of the matrix when the width of the network layers is large. Our assumptions on the dimensions of matrices $W_\ell$ will be phrased in terms of their relation to the sample size $n$, i.e., the number of columns of the \emph{data} matrix $X$. The primary range of interest is the so called random matrix regime, when $n_0,n_1,\ldots,n_L$ are each proportional to $n$, however our inequalities will be also applicable in the setting when for instance $n_0$, i.e., the dimension of the data, is fixed and the sample size is large, or the data is high-dimensional ($n_0$ is large) and for instance $n = 1$, which corresponds to the action of the network on a single vector.

Let us also mention that we consider a situation when the data matrix $X$ is random and well-concentrated. Many results in neural network literature consider instead a situation when $X$ is deterministic, but satisfies some normalization property, such as restrictions on the operator norm. Under appropriate normalization of the operator norm of the data matrix $X$ (of order roughly $\sqrt{n_0} + \sqrt{n}$) our arguments give results also in such models, by a very mild change of the proofs. However, since already for random $X$, our results are split into several cases, depending on the type of assumptions on the network, and as explained, mathematically the case of deterministic $X$ does not involve any new ingredients, we do not pursue this direction.

\subsection{Bounded activation function}\label{sec:bounded-activation}

In this section we will be working under the following set of assumptions.

Recall that we consider $\XX = (X,(W_\ell,B_\ell)_{\ell=0,\ldots,L-1})$, where $X$ is an $n_0 \times n$ random matrix, $W_\ell$ is an $n_{\ell+1}\times n_\ell$ random matrix and $B_\ell$ is an $n_{\ell+1}$-dimensional random vector. When discussing the concentration properties of $\XX$, we view it as a random vector with values in $\R^{N}$, where $N = nn_0 + \sum_{\ell=0}^{L-1} (n_\ell n_{\ell+1} + n_{\ell+1})$. The random matrix $Z_L$ is defined via the recursive formula \eqref{eq:model}.

Our simplest result will be derived under the following set of assumptions.

\textbf{Assumption A1}
\begin{enumerate}
  \item $\XX$ satisfies the concentration property \eqref{eq:concentration-assumption}.
  \item $X$ and $W_\ell$, $\ell=0,\ldots,L-1$ are centered.
  \item For all $\ell = 0,\ldots, L-1$, the activation function $\sigma_\ell$ is $\lambda_\sigma$-Lipschitz and $\|\sigma_\ell\|_\infty \le \kappa$.
  \item  For all $\ell = 1,\ldots,L$, $n/n_\ell \le \alpha$.
  \item For all $\ell = 1,\ldots,L$, $n_{\ell-1}/n_{\ell} \le \beta$.
  \item The data dimension $n_0$ satisfies $n_0 \ge \eta$.
\end{enumerate}

The hypotheses in Assumption A1 can be divided into three sets: assumptions on the distribution of $\XX$ (items (1),(2)), assumptions on the activation functions (item (3)), and assumptions on the dimensions of the matrices, i.e., the sample size $n$, the dimension of the data $n_0$ and the width of consecutive layers of the network $n_1,\ldots,n_L$ (items (4)--(6)). In subsequent sections we will consider various modifications of Assumption A1, each time their structure will be similar. Let us note that trivially one can always take $\eta = 1$, while the parameters $\alpha,\beta$ may take arbitrary positive values. Moreover, the assumption A1 becomes stronger with increasing $\eta$ and decreasing $\alpha$ and $\beta$. From the machine learning perspective, one can say that $L,\lambda_\sigma$, and $\kappa$ depend on the general \emph{architecture} of the network, while $n_\ell$, $\ell = 1,\ldots,L$ are the widths of its layers, chosen to be sufficiently large to allow for a small (in the overparametrized regime even null) empirical loss while training on $n$ samples of dimension $n_0$.

\begin{theorem}\label{thm:bounded-activation}
There exists a constant $c>0$, depending only on $L,\rho,\lambda_\sigma,\kappa,\alpha,\beta,\eta$, such that if the Assumption A1 holds, then for every function $F \colon \R^{n_L \times n} \to \R$, $1$-Lipschitz with respect to the Hilbert--Schmidt norm, and every $t > 0$,
\begin{align}\label{eq:bounded-activation-concentration-simplified-subgaussian}
 \p(|F(Z_L) - \Med F(Z_L)| \ge t) \le 2\exp(-ct^2).
\end{align}
Moreover if by $c_{opt} = c_{opt}(L,\rho,\lambda_\sigma,\kappa,\alpha,\beta,\eta)$ we denote the optimal (the largest) constant in \eqref{eq:bounded-activation-concentration-simplified-subgaussian}, then if simultaneously $\alpha\to 0$ and $\eta \to \infty$, while the other parameters are fixed, we have $c_{opt} \to \infty$.

\end{theorem}

\begin{remark}
At the cost of altering the constant $c$ by a universal factor, one can replace the median in \eqref{eq:bounded-activation-concentration-simplified-subgaussian} by the expectation $\E F(Z_L)$ (see Lemma \ref{le:concentration-around-mean}).
\end{remark}

\begin{remark}
Examples of activation functions considered in the neural network literature and covered by Assumption 1 are, e.g., the sigmoid function, $\sin$, $\tanh$.
\end{remark}
\begin{remark}An important aspect of the above theorem is that it is dimension-free, in the sense that for fixed values of $L,\rho,\lambda_\sigma,\kappa,\alpha,\beta,\eta$, the right-hand side of the estimate does not depend on $n$. In considerations related to random matrices one usually assumes that $n_\ell/n \to \alpha_\ell$ for some positive constants $\alpha_\ell$, which makes our results applicable in this setting. In Sections \ref{sec:symmetry} under additional assumptions we provide inequalities which improve as the widths of the network tend to $\infty$ while the parameter $\alpha$ is fixed. In Section \ref{sec:discussion} we show that under our general assumptions dimension-free concentration is the best one that can be achieved. We also discuss the dependence of our estimates on parameters $\alpha,\beta$ and we show that if $\alpha \to 0$ but $\eta \nrightarrow \infty$, then the constant $c$ remains bounded from above.
\end{remark}

Let us conclude this section with an application to linear eigenvalue statistics of the random matrix $Z_L$. It is a direct consequence of Theorem \ref{thm:bounded-activation} and the Lipschitz property of such functions, mentioned below formula \eqref{eq:les}.

\begin{cor}\label{cor:lse-bounded-activation} In the setting of Theorem \ref{thm:bounded-activation}, if $f\colon \R \to \R$ is 1-Lipschitz and $S_f(Z_L)$ is given by \eqref{eq:les}, then for every $t > 0$,
\begin{displaymath}
  \p\Big(|S_f(Z_L) - \E S_f(Z_L)|\ge t\Big) \le 2\exp(-c t^2/n).
\end{displaymath}
\end{cor}

\subsection{General Lipschitz activation function}\label{sec:Lipschitz-activation}

We will now generalize Theorem \ref{thm:bounded-activation} to the setting of unbounded activation functions. More precisely, let us consider the following set of assumptions.

\textbf{Assumption A2}
\begin{enumerate}
  \item $\XX$ satisfies the concentration property \eqref{eq:concentration-assumption}.
  \item $X$ and $W_\ell, B_\ell$, $\ell=0,\ldots,L-1$, are centered.
  \item For all $\ell = 0,\ldots, L-1$, the activation function $\sigma_\ell$ is $\lambda_\sigma$-Lipschitz and $|\sigma_\ell(0)|\le \tau$.
  \item For all $\ell = 1,\ldots,L$, $n/n_\ell \le \alpha$.
  \item For all $\ell = 1,\ldots,L$, $n_{\ell-1}/n_{\ell} \le \beta$.
  \item The data dimension $n_0$ satisfies $n_0 \ge \eta$.

\end{enumerate}

\begin{remark}
The only difference between the Assumptions A1 and A2 is the item (3), in Assumption A2 we allow the activation functions to be unbounded. In particular this assumption is satisfied by the ReLU function or the absolute value, which are widely used in practice.
\end{remark}

\begin{theorem}\label{thm:unbounded-activation}
 There exists a constant $c>0$, depending only on $L,\rho,\lambda_\sigma,\tau,\alpha,\beta,\eta$, such that if the Assumption A2 holds, then for every function $F \colon \R^{n_L \times n} \to \R$, $1$-Lipschitz with respect to the Hilbert--Schmidt norm, and every $t > 0$,
\begin{align}\label{eq:unbounded-activation-concentration}
 \p(|F(Z_L) - \Med F(Z_L)| \ge t) \le 2\exp\Big(-c\min_{1\le k \le L+1} \Big( n^{(k-1)/k}t^{2/k}\Big)\Big).
\end{align}
If in addition $\|\sigma_{L-1}\|_\infty \le \kappa$, then
\begin{align}\label{eq:unbounded-but-the-last-one}
  \p(|F(Z_L) - \Med F(Z_L)| \ge t) \le 2\exp(-c' t^2)
\end{align}
for some positive constant $c'$, depending only on $L,\rho,\lambda_\sigma,\tau,\kappa,\alpha,\beta,\eta$. Moreover, if simultaneously $\alpha \to 0$ and $\eta \to \infty$, while the other parameters are fixed, the optimal constants $c_{opt},c'_{opt}$ in the above inequalities diverge to $\infty$.
\end{theorem}

\begin{remark} An analogous result holds for concentration around the mean (see Lemmas \ref{le:concentration-around-mean} and \ref{le:concentration-equivalence}).
\end{remark}

We also have the following proposition, which will in fact be proven inductively alongside with Theorem \ref{thm:unbounded-activation}.

\begin{prop}\label{prop:unbounded-activation} In the setting of Theorem \ref{thm:unbounded-activation}, there exists a constant $C$, depending only on $L,\rho,\lambda_\sigma,\tau,\alpha,\beta,\eta$, such that
\begin{displaymath}
  \E\|Z_L\|_{HS} \le C\sqrt{n}.
\end{displaymath}
\end{prop}

In the language of neural networks, in Theorems \ref{thm:bounded-activation} and \ref{thm:unbounded-activation} we require that the width of consecutive layers be at least of the order of the sample size (as described by the coefficient $\alpha$) and that each consecutive layer of the network be of width at least of the order of the previous one (which is captured by the coefficient $\beta$). If we impose additional independence conditions, which are typical for neural networks, one can eliminate the coefficient $\beta$. The corresponding set of assumptions is

\textbf{Assumption A3}
\begin{enumerate}
  \item The random elements $X, (W_\ell,B_\ell)$, $\ell=0,\ldots,L-1$ are stochastically independent.
  \item Each of the elements $X, (W_\ell,B_\ell)$, $\ell=0,\ldots,L-1$ (seen as a random vector in $\R^{n_0n}$ and $\R^{n_{\ell+1}n_\ell + n_{\ell+1}}$ respectively) satisfies the concentration property \eqref{eq:concentration-assumption}.
  \item $X$ and $W_\ell, B_\ell$, $\ell=0,\ldots,L-1$, are centered.
  \item For all $\ell = 0,\ldots, L-1$, the activation function $\sigma_\ell$ is $\lambda_\sigma$-Lipschitz and $|\sigma_\ell(0)|\le \tau$.
  \item For all $\ell = 1,\ldots,L$, $n/n_\ell \le \alpha$.
  \item The data dimension $n_0$ satisfies $n_0 \ge \eta$.
\end{enumerate}

Thus in Assumption A3 one imposes stronger conditions on the distribution of $\XX$ (items (1)--(3)) than in Assumption A2 (recall Proposition \ref{prop:independent-concentration}, which implies that items (1) and (2) of Assumption A3 imply item (1) of Assumption A2 with an appropriate change of the parameter, depending only on $L$), but weakens the requirements on the widths of the layers by eliminating the coefficient $\beta$.
\begin{theorem}\label{thm:unbounded-activation-independence}
  There exists a constant $c$, depending only on $L,\rho, \lambda_\sigma,\tau, \alpha, \eta$, such that whenever the Assumption A3 holds, then for every $F\colon \R^{n_L\times n} \to \R$, 1-Lipschitz with respect to the Hilbert--Schmidt norm, and every $t > 0$,
\begin{align}\label{eq:unbounded-activation-concentration-independence}
 \p(|F(Z_L) - \E F(Z_L)| \ge t) \le 2\exp\Big(-c\min_{1\le k \le L+1} \Big( n^{(k-1)/k}t^{2/k}\Big)\Big).
\end{align}
If in addition $\|\sigma_{L-1}\|_\infty \le \kappa$, then
\begin{align}\label{eq:unbounded-but-the-last-one-independence}
  \p(|F(Z_L) - \E F(Z_L)| \ge t) \le 2\exp(-c' t^2)
\end{align}
for some positive constant $c'$, depending only on $L,\rho,\lambda_\sigma,\tau,\kappa,\alpha,\eta$. Moreover, if simultaneously $\alpha \to 0$ and $\eta \to \infty$, while the other parameters are fixed, the optimal constants $c_{opt},c'_{opt}$ in the above inequalities diverge to $\infty$.\end{theorem}

\begin{remark} An analogous result holds for concentration around the mean (see Lemmas \ref{le:concentration-around-mean} and \ref{le:concentration-equivalence}).
\end{remark}

Again, we also have an estimate on the expected Hilbert--Schmidt norm of $Z_L$.

\begin{prop}\label{prop:unbounded-activation-independence} In the setting of Theorem \ref{thm:unbounded-activation-independence}, there exists a constant $C$, depending only on $L,\rho, \lambda_\sigma,\tau, \alpha, \eta$, such that
\begin{displaymath}
  \E\|Z_L\|_{HS} \le C\sqrt{n}.
\end{displaymath}
\end{prop}

\begin{cor}\label{cor:les-unbounded-activation} In the setting of Theorem \ref{thm:unbounded-activation}, if $f\colon \R \to \R$ is 1-Lipschitz and $S_f(Z_L)$ is given by \eqref{eq:les}, then for every $t > 0$,
\begin{displaymath}
  \p\Big(|S_f(Z_L) - \E S_f(Z_L)|\ge t\Big)  \le 2\exp\Big(-c\Big(\min_{1\le k \le L+1} n^{(k-2)/k} t^{2/k}\Big) \Big)
\end{displaymath}
If in addition $\|\sigma_{L-1}\|_\infty \le \kappa$, then
\begin{displaymath}
    \p\Big(|S_f (Z_L)- \E S_f(Z_L)|\ge t\Big) \le 2\exp(-c't^2/n).
\end{displaymath}
An analogous statement holds in the setting of Theorem \ref{thm:unbounded-activation-independence}.
\end{cor}

\subsection{Additional symmetry}\label{sec:symmetry}

The inequalities presented in the previous sections in the situation when all $n_\ell$'s are proportional to $n$, provide concentration whose dependence on $n,n_0$ does not improve substantially when $n,n_0\to \infty$, i.e., for each fixed $t$, eventually the tail estimate becomes of the form $2\exp(-ct^2)$, which is separated from zero. This is in contrast with the inequalities for Wigner matrices and is related to the fact that the matrices $Z_\ell$, which act \emph{multiplicatively} on the concentrated matrices $W_\ell$, have operator norm of order $\sqrt{n}$. Similarly, the deterministic vectors $\1_n$ which multiply the concentrated vectors $B_n$ have Euclidean norm equal to $\sqrt{n}$. Below we consider a set of assumption, under which this phenomenon does not occur. In particular we consider a situation without bias and with additional symmetry of the distribution of $W_\ell$ and the activation function, which will allow to show that $\|Z_\ell\|_{op} = \mathcal{O}(1)$, leading to concentration rates which improve with $n,n_0\to \infty$. This counterparts the results obtained in the classical Wigner setting by Guionnet and Zeitouni \cite{MR1781846}.

\medskip

\textbf{Assumption A4}
\begin{enumerate}
  \item The random elements $X, W_\ell$, $\ell=0,\ldots,L-1$ are stochastically independent.
  \item Each of the elements $X, W_\ell$, $\ell=0,\ldots,L-1$ (seen as a random vector in $\R^{n_0n }$ and $\R^{n_{\ell+1} n_\ell}$ respectively) satisfies the concentration property \eqref{eq:concentration-assumption}.
  \item $X$ is centered.
  \item For all $\ell = 0,\ldots, L-1$, $B_\ell = 0$.
  \item For all $\ell=0,1,\ldots,L-1$, the law of $W_\ell$ is symmetric with respect to the origin.
  \item For all $\ell = 0,\ldots, L-1$, the activation function $\sigma_\ell$ is $\lambda_\sigma$-Lipschitz and odd,
  \item For all $\ell = 1,\ldots,L$, $n/n_\ell \le \alpha$.
\end{enumerate}

\begin{theorem}\label{thm:unbounded-activation-symmetry-independence}
There exists a constant $c$, depending only on $L,\rho,\lambda_\sigma, \alpha$, such that whenever the assumption A4 holds, then for every $F\colon \R^{n_L\times n}\to \R$, $1$-Lipschitz with respect to the Hilbert--Schmidt norm, and every $t>0$,
\begin{align}\label{eq:unbounded-activation-symmetry-independence}
\p(|F(Z_L) - \E F(Z_L)| \ge t) \le 2\exp\Big(-c\min(n,n_0)\min(t^2,t^{2/(L+1)})\Big).
\end{align}
\end{theorem}

\begin{remark} An analogous result holds for concentration around the median (see Lemma \ref{le:concentration-equivalence}).
\end{remark}

\begin{remark}\label{re:no-n}
An easy modification of the proof of the above theorem shows that under an additional assumption that $n_0/n_\ell \le \gamma$ for $\ell=1,\ldots,L$, one can replace $\min(n,n_0)$ simply by $n_0$ at the cost of making the constant $c$ dependent additionally also on $\gamma$. Clearly without this assumption we cannot eliminate $n$ from the inequality, since for $L = 1$, $n = n_1 = 1$, if we choose the entries of $W_0$ and $X$ to be i.i.d. standard Gaussian variables and set $F(x) = \sigma_0(x) = x$, we get
\begin{displaymath}
  \Var(F(Z_1)) = \Var\Big(\sum_{i=1}^{n_0} \frac{1}{\sqrt{n_0}} W(1,i)X(i,1)\Big) = 1,
\end{displaymath}
which is independent of $n_0$.
\end{remark}

The following proposition gives an estimate on the expected operator norm of the matrix $Z_L$.

\begin{prop}\label{prop:unbounded-activation-symmetry-independence} In the setting of Theorem \ref{thm:unbounded-activation-symmetry-independence}, there exists a constant $C$, depending only on $L,\rho, \lambda_\sigma,\alpha$, such that
\begin{displaymath}
  \E\|Z_L\|_{op} \le C\Big(1 + \sqrt{\frac{n}{n_0}}\Big).
\end{displaymath}
\end{prop}

Let us also state improved (with respect to Corollaries \ref{cor:lse-bounded-activation} and \ref{cor:les-unbounded-activation}) concentration estimates for linear eigenvalue statistics.

\begin{cor}\label{cor:lse-unbounded-activation-independence-symmetry} In the setting of Theorem \ref{thm:unbounded-activation}, if $f$ is 1-Lipschitz and $S_f(Z_L)$ is given by \eqref{eq:les}, then for every $t > 0$,
\begin{displaymath}
  \p\Big(|S_f(Z_L) - \E S_f(Z_L)|\ge t\Big)  \le 2\exp\Big(-c  \min(n,n_0) \min\Bigl(\frac{t^2}{n}, \frac{t^{2/(L+1)}}{n^{-1/(L+1)}}\Bigr)\Big).
\end{displaymath}
\end{cor}

In the random matrix regime, when all $n_\ell$'s $\ell=0,\ldots,L$ are proportional to $n \to \infty$, the above corollary gives
\begin{displaymath}
  \p\Big(|S_f(Z_L) - \E S_f(Z_L)|\ge t\Big) \le 2\exp(-\tilde{c} \min(t^2,n^{L/(L+1)} t^{2/(L+1)}).
\end{displaymath}
In particular, by integration by parts, $\Var(S_f(Z_L)) = \mathcal{O}(1)$ as $n\to \infty$. This suggest that under additional assumptions on the entries (e.g., independence), $S_f(Z_L)$ should satisfy a CLT without any additional normalization, as it is the case for classical random matrix ensembles (see, e.g., \cite{MR2561434,MR2040792,MR2567175,MR2808038}).

\subsection{Non-commutative polynomials in random matrices}\label{sec:polynomials}

In this section we will show how arguments similar to those presented above can be used to obtain concentration for non-commutative polynomials in concentrated random matrices. Since this model is different from the one considered in previous sections, we will change the notation and denote the random element consisting of all the random matrices involved by $\xX$ and not by $\XX$. We will also allow complex entries. Below, by $A^\ast$ we denote the conjugate of a matrix $A$, i.e., $A^\ast = \bar{A}^T$.

\begin{theorem}\label{thm:polynomials}
  Let $\xX = (X_1,\ldots,X_L)$, where $X_i$ are $n\times n$ centered complex random matrices. Assume that $\xX$ satisfies the subgaussian concentration inequality \eqref{eq:concentration-assumption}. Let $P\colon (\C^{n\times n})^L\to \C$ be a non-commutative polynomial of degree $d$, of the form
  \begin{align}\label{eq:polynomial}
    P(x_1,\ldots,x_L) = \sum_{k=0}^d \sum_{\ii = (i_1,\ldots,i_k) \in [L]^k} \sum_{\varepsilon = (\varepsilon_1,\ldots,\varepsilon_k) \in \{1,\ast\}^k} A_{k,\ii,\varepsilon}^{(0)}x_{i_1}^{\varepsilon_1}A_{k,\ii,\varepsilon}^{(1)} x_{i_2}^{\varepsilon_2}\cdots A_{k,\ii,\varepsilon}^{(k-1)}x_{i_k}^{\varepsilon_k} A_{k,\ii,\varepsilon}^{(k)},
  \end{align}
  where $A_{k,\ii,\varepsilon}^{(j)}$ are $n\times n$ deterministic matrices with operator norm bounded by 1. Let $\yY = P(n^{-1/2}\xX)$. Then, for any function $F \colon \C^{n\times n} \to \R$, $1$-Lipschitz with respect to the Hilbert--Schmidt norm, and any $t > 0$,
\begin{displaymath}
  \p(|F(\yY) - \Med F(\yY)| \ge t) \le 2\exp\Big(- c n \min(t^2,t^{2/d})\Big),
\end{displaymath}
where $c$ is a constant depending only on $d,L,\rho$.
\end{theorem}

\begin{remark}
The assumptions of the above theorem are satisfied for instance if the matrices $X_1,\ldots,X_N$ are independent and each of them satisfies the subgaussian concentration property with constant $\rho'$. In this case one can take $\rho = C_L \rho'$ (see Proposition \ref{prop:independent-concentration}). Another model for which the assumptions are satisfied is given by multimatrix models defined by \eqref{eq:multi}. We refer to \cite{MR3646880, MR2353386} for convexity assumptions on the potential $V$ and the cut-off level $M$, which guarantee the concentration property.
\end{remark}

\begin{cor} In the setting of Theorem \ref{thm:polynomials}, assume that the polynomial $P$ is formally self-adjoint. Let $\gamma_1,\ldots,\gamma_n$ be the eigenvalues of $\yY$ and for $f\colon \R \to \R$ define
\begin{displaymath}
  S_f = \sum_{i=1}^n f(\gamma_i).
\end{displaymath}
If $f$ is 1-Lipschitz, then for every $t\ge 0$,
\begin{align}\label{eq:les-polynomials}
  \p(|S_f- \E S_f| \ge t) \le 2\exp\Big(- c  \min(t^2,n^{1-1/d} t^{2/d})\Big).
\end{align}
\end{cor}

\begin{remark}\label{le:polynomials}
If $f(x) = x$, i.e., $S_f = \tr \yY$, then an easy modification of the proof of Theorem \ref{thm:polynomials} allows to replace $n^{1-1/d}$ in the above corollary by $n$. The assumption that $P$ is formally self-adjoint is not needed.

Such an inequality for the trace was obtained by Meckes and Szarek \cite{MR2869165} under a weaker assumption of convex concentration property (i.e., concentration property for \emph{convex} Lipschitz functions). In the case of Wigner matrices (i.e., Hermitian matrices with independent entries on and above the main diagonal) with subexponential entries related inequalities for the trace have been very recently obtained by Bao and Munoz George \cite{bao2024ultrahighordercumulants} under the assumption that the matrices $X_i$ are independent. We are not aware of inequalities for general linear spectral statistics of polynomials in random matrices.

The tail estimate in the result by Meckes and Szarek is of the form $\exp(-c\min(t^2,n t^{2/d})$, which has a better dependence on $n$ than in our inequality \eqref{eq:les-polynomials}. Thus in this special case of $f(x) = x$ our proof allows to recover this estimate, see Remark \ref{re:trace} following the proof of Theorem \ref{thm:polynomials} in Section \ref{sec:polynomials-proofs} (however, as mentioned above, our assumptions are stronger than those from \cite{MR2869165}). Bao and Munoz George have the tail of the form $\exp(-c\frac{t^2}{1 + t^{2-2/d}n^{-2/d}})$ in the case of GUE/GOE, which is weaker than the one in the result by Meckes--Szarek and in fact for  $d\ge 4$ weaker than our estimate \eqref{eq:les-polynomials}. They also provide a different type of bound for the trace of Wigner matrices, satisfying certain assumptions on the cumulants.
\end{remark}

\subsection{Further applications}\label{sec:applications}

Let us now formulate some additional corollaries to our main inequalities. Their proofs will be postponed to Section \ref{sec:applications-proofs}.

The first corollary concerns convergence of the empirical spectral distribution of the conjugate kernel $Z_L^TZ_L$. Recall that the spectral measure of an $n\times n$ symmetric matrix $A$ is the Borel probability measure on $\R$, defined as
\begin{displaymath}
  \mu_A = \frac{1}{n}\sum_{i=1}^n \delta_{\gamma_i(A)},
\end{displaymath}
where $\gamma_i(A)$ $i = 1,\ldots,n$, are the eigenvalues of $A$, counted with multiplicities and $\delta_x$ stands for the Dirac mass at $x$. If $A$ is a random matrix, then $\mu_A$ is a random probability measure.

From Theorem \ref{thm:unbounded-activation} we obtain the following corollary.
\begin{cor}\label{cor:distance-bound} Fix a nonnegative integer $L$. Consider a sequence of positive integers $(n^{(m)})_{m\in\N}$, such that $\lim_{m\to \infty} n^{(m)} = \infty$, a corresponding sequence of $(L+1)$-tuples of positive integers $(n_\ell^{(m)})_{\ell = 0}^{L}$, Lipschitz activation functions $(\sigma_{\ell}^{(m)})_{\ell=0}^{L-1}$, and a sequence of random elements $\XX^{(m)} = (X^{(m)},(W^{(m)}_\ell,B^{(m)}_{\ell})_{\ell=0}^{L-1})$, where $X^{(m)}$ is an $n^{(m)}_0\times n^{(m)}$ random matrix, $W^{(m)}_\ell$ is an $n_{\ell}^{(m)}\times n_{\ell+1}^{(m)}$ random matrix and $B_{\ell}^{(m)}$ is an $\R^{n_{\ell+1}^{(m)}}$-valued random vector.

Assume that for each $m$, $\XX^{(m)}$ satisfies the Assumption 2 (resp. the Assumption 3) with the parameters $\rho, \lambda_\sigma,\tau,\alpha,\beta$ (resp. $\rho,\lambda_\sigma,\tau,\alpha$) independent of $m$. Finally let $A_m = (Z_L^{(m)})^T Z_L^{(m)}$ be the corresponding conjugate kernel. Then for any distance $d$ metrizing the weak convergence of probability measures on $\R$, $d(\mu_{A_m},\E \mu_{A_m})$ converges to zero in probability as $m\to \infty$.
\end{cor}

Another corollary concerns concentration for linear eigenvalue statistics of the conjugate kernel $A = Z_L^T Z_L$. In Corollaries \ref{cor:lse-bounded-activation}, \ref{cor:les-unbounded-activation} and \ref{cor:lse-unbounded-activation-independence-symmetry} we considered linear spectral statistics $S_f(Z_L)$, where $f\colon \R \to \R$ is a Lipschitz function (recall the definition \eqref{eq:les}). Clearly $S_f(A)$ can be written as $S_g(Z_L)$ for $g(x) = f(x^2)$. However the function $g$ in general is not Lipschitz, therefore we need one additional step to pass to $S_g$. One can use to this end a result from \cite{louart2025operationconcentrationinequalities}, however, since our setting is simpler than the general one presented in this article, and moreover using \cite{louart2025operationconcentrationinequalities} would require additional involved notation, we will provide a direct proof. For simplicity we will state the result only in the case where the layers are independent with widths (including the data dimension) proportional to $n$.

\begin{cor}\label{cor:les-conjugate-kernel}
Assume that the Assumption A4 holds and that additionally $n/n_0 \le \alpha$. Let $Y$ be the $n$-dimensional random vector whose components are the eigenvalues of the conjugate kernel $A=Z_L^TZ_L$. Then for any $F \colon \R^n \to \R$, 1-Lipschitz with respect to the standard Euclidean norm, and any $t > 0$,
\begin{displaymath}
  \p(|F(Y) - \E F(Y)| \ge t) \le 2\exp\Big( - c n\min(t^2, t^{1/(L+1)})\Big),
\end{displaymath}
where $c$ is a constant depending on $L,\rho,\lambda_\sigma,\alpha$.

In particular,
\begin{displaymath}
  \p(|S_f(A) - \E S_f(A)| \ge t) \le 2\exp\Big(-c\min(t^2,n^{(2L+1)/(2L+2)}t^{1/(L+1)})\Big).
\end{displaymath}
\end{cor}

\begin{remark} It is easy to see that the exponent $t^{1/(L+1)}$ in the above inequality cannot be improved, since in the case $\sigma_\ell(x) = x$ if all the involved matrices are independent with standard Gaussian entries, the entries of the conjugate kernel are polynomials of degree $2(L+1)$ of independent standard Gaussian variables, and for Gaussian polynomials the tail for large values of $t$ behaves as $\exp(-c t^{2/d})$, where $d$ is the degree (which agrees with Theorem \ref{thm:polynomials}).

Again, the fact that the above inequality yields an $\mathcal{O}(1)$ bound on the variance suggests that the linear eigenvalue statistics of the conjugate kernel may satisfy a central limit theorem (perhaps under additional smoothness assumptions on the activation functions).
\end{remark}

Our final application is a deviation inequality for the Wasserstein distance $\mathcal{W}_1$ between the expected and empirical spectral distribution of the conjugate kernel  (recall that definition \eqref{eq:Wasserstein} of the Wasserstein distance).  By Kantorovich duality (see, e.g., \cite{MR1964483}), one can equivalently express $\mathcal{W}_1(\mu,\nu)$ as
\begin{align}\label{eq:Wasserstein-duality}
  \mathcal{W}_1(\mu,\nu) = \sup_{f\colon \R \to \R, \|f\|_{Lip} \le 1} \Big|\int_{\R} f d\mu - \int_\R fd\nu\Big|,
\end{align}
with supremum running over all $1$-Lipschitz functions $f\colon \R\to \R$.

The following Proposition is inspired by results from \cite{MR3055255}, where the authors discuss convergence rates in the Wasserstein distance for many classical random matrix models satisfying subgaussian concentration. The proof, presented in Section \ref{sec:applications-proofs} is an adaptation of their arguments.

\begin{prop}\label{prop:Wasserstein-convergence}
Fix $r \in (0,2]$ and let $A$ be an $n \times n$ random Hermitian matrix, such that $\E \|A\|_{op} \le D$ and for all $F \colon \R^{n\times n} \to \R$, and all $t > 0$,
\begin{align}\label{eq:Wasserstein-concentration-assumption}
  \p(|F(A) - \E F(A)| \ge t) \le 2\exp(-an\min(t^2,t^r)).
\end{align}
Then $\E \mathcal{W}_1(\mu_A,\E \mu_A) \le C/n^{2/3}$ and for all $t > 0$,
\begin{align*}
  \p\Big(\mathcal{W}_1(\mu_A,\E \mu_A) \ge C/n^{2/3} + t\Big)\le 2\exp(-a\min(n^2t^2,n^{1+r/2}t^{r})),
\end{align*}
where $C$ depends only on $a,D,r$.
\end{prop}

Applying the above proposition  to the conjugate kernel, we obtain
\begin{cor}\label{cor:Wasserstein-conjugatekernel}
In the setting of Corollary \ref{cor:les-conjugate-kernel}, $\E \mathcal{W}_1(\mu_A,\E \mu_A) \le C/n^{2/3}$ and for all $t > 0$,
\begin{displaymath}
  \p(\mathcal{W}_1(\mu_A,\E \mu_A) \ge C/n^{2/3} + t)\le 2\exp(-c \min(n^2 t^2, n^{\frac{2L+3}{2L+2}}t^{\frac{1}{L+1}})),
\end{displaymath}
where $C,c$ are positive constants, depending only on $L,\rho,\lambda_\sigma,\alpha$.
\end{cor}

We can also apply Proposition \ref{prop:Wasserstein-convergence} to non-commutative polynomials, obtaining the following result.

\begin{cor}\label{cor:Wasserstein-polynomials}In the setting of Theorem \ref{thm:polynomials}, assume that the polynomial $P$ is formally self-adjoint. Then,
$\E \mathcal{W}_1(\mu_\yY,\E \mu_\yY) \le C/n^{2/3}$ and for all $t > 0$,
\begin{displaymath}
  \p(\mathcal{W}_1(\mu_\yY,\E \mu_\yY) \ge C/n^{2/3} + t) \le 2\exp(-c\min(n^2 t^2,n^{\frac{d+1}{d}}t^{2/d})),
\end{displaymath}
where $C,c$ are positive constants depending only on $d,L,\rho$.
\end{cor}

\subsection{Discussion of scaling and dependence on parameters}\label{sec:discussion}

As already mentioned, the right-hand side of inequalities presented in Theorems \ref{thm:bounded-activation}, \ref{thm:unbounded-activation}, \ref{thm:unbounded-activation-independence} does not tend to zero with the dimensions of the data tending to infinity, which is in contrast with the estimates for classical random matrices with independent Gaussian entries of variance $1/n$. Below we will argue that for general Lipschitz activation functions one cannot hope for such an improvement, even if the activation functions are bounded, there is no bias, and the entries of all the matrices involved are i.i.d. random variables. We will also discuss the role of the parameters relating the dimensions of the matrices appearing in the definition of $Z_L$.

\subsubsection{Dependence on $\alpha$} We will start by discussing the dependence of our estimates on the parameter $\alpha$, which controls the ratio $n/n_\ell$ for $\ell = 1,\ldots,L$, i.e., the relation between the sample size and the minimum width of the network's layer.

\subsubsection*{Networks with bias} In the presence of bias, the dependence of the constants in all our inequalities on the parameter $\alpha$ is clearly necessary. Even in the case when $L=1$, $\sigma_0(x) = \min(\max(x,-1),1)$, $W_0,X = 0$ and $B_0$ has independent coordinates distributed uniformly on $(-1,1)$, for the function $F\colon \R^{n\times n_1} \to \R$ given by $F(z) = \frac{1}{\sqrt{nn_1}}\sum_{i=1}^{n_1}\sum_{j=1}^n z_{ij}$ (which is $1$-Lipschitz), we have
\begin{displaymath}
  F(Z_1) = \frac{\sqrt{n}}{n_1} \sum_{i=1}^{n_1} B_0(i)
\end{displaymath}
and so $\Var(F(Z_1)) = \frac{n}{3n_1}$. Thus, if $n/n_1 \to \infty$ (i.e., the coefficient $\alpha$ becomes unbounded), the variance diverges and so, $F(Z_1)$ does not satisfy concentration with a dimension-independent profile. For $n/n_1$ fixed the variance is also of fixed order and thus the concentration rate cannot improve with $n$.

\subsubsection*{No bias, growing $\eta$} Let us now argue that the dependence on $\alpha$ is needed even if there is no bias and the matrices $X,W_\ell$ are independent. We will need to consider $L \ge 2$ since for $L=1$, one can take advantage of a better control over $\|Z_0\|_{op}$, which in the absence of bias leads to improved estimates (see Remark \ref{re:L1}).

Consider thus $L = 2$, $n_0 = n_1\to \infty $ as $n\to \infty$ and $n_2 = an$. Let $\sigma_0(x) = \sigma_1(x) = \sigma(x) := \min(1,\max(x,0))$ (so that the parameters $\lambda_\sigma, \kappa$ of the Assumptions A1, A2, A3 satisfy $\lambda_\sigma = \kappa = 1$), $F(z) = \frac{1}{\sqrt{n}}\sum_{j=1}^n z_{1,j}$ (note that $F$ is 1-Lipschitz) and assume that the entries of $X,W_0,W_1$ are i.i.d. standard Gaussian variables. Set also $B_\ell=0$ for all $\ell$.

It turns out  that for large $n$,
\begin{align}\label{eq:variance-behaviour}
\Var(F(Z_L)) \ge \frac{c}{a},
\end{align}
where $c$ is a positive constant. Note that if $n_0 = n_1 \le a n$, then the parameter $\beta$
in Assumptions A1, A2 is bounded from above by 1. Also the parameters concerning the
activation function and the concentration parameter $\rho$ are bounded. On the other hand, the parameter $\alpha$ diverges when $a \to 0$.

This shows that the constants $c$ in inequalities \eqref{eq:bounded-activation-concentration-simplified-subgaussian} of Theorem \ref{thm:bounded-activation}, \eqref{eq:unbounded-activation-concentration} of Theorem \ref{thm:unbounded-activation}, and \eqref{eq:unbounded-activation-concentration-independence} of Theorem \ref{thm:unbounded-activation-independence} must depend on $\alpha$, even in the simplest case, when there is no bias and all the matrices are independent, with i.i.d.
standard Gaussian entries.

Below we provide a sketch of the proof of \eqref{eq:variance-behaviour}, leaving certain details to the reader.

We have
  \begin{displaymath}
    Z_2(i,j) = \frac{1}{\sqrt{an}}\sigma\Big(\frac{1}{\sqrt{n_0}} \sum_{k_1=1}^{n_0} W_1(i,k_1)\sigma\Big(\frac{1}{\sqrt{n_0}}\sum_{k_0=1}^{n_0} W_0(k_1,k_0)X(k_0,j)\Big)\Big).
  \end{displaymath}
We will first show that for any $(i_1,j_1) \neq (i_2, j_2)$, as $n \to \infty$, the random vector
 \begin{align}\label{eq:convergence-in-law}
   \sqrt{an}(Z_2(i_1,j_1),Z_2(i_2,j_2)) \textrm{ converges in distribution to } (\sigma(G_1),\sigma(G_2)),
 \end{align}
where $(G_1,G_2)$ is a mean-zero Gaussian vector in $\R^2$, with $\E G_1^2 = \E G_2^2 = \E\sigma^2(g)$, $\E G_1G_2 = \ind{i_1=i_2} (\E \sigma(g))^2$ for a standard Gaussian variable $g$. In particular $G_1,G_2$ are independent for $i_1\neq i_2$ and are positively correlated if $i_1=i_2$. In the following discussion we allow $a$ to be arbitrary (possibly depending on $n$), we just assume that $n_0=n_1\to \infty$, $n_2 = an$.

To justify the convergence, observe that conditionally on $X$ and $W_0$, the random vector $\sqrt{an}(Z_2(i_1,j_1),Z_2(i_2,j_2))$ is distributed as $(\sigma(G_1^{(n)}),\sigma(G_2^{(n)}))$, where $(G_1^{(n)},G_2^{(n)})$ is a centered Gaussian vector such that
 \begin{align}
   \label{eq:first} \E (G_m^{(n)})^2 & = \frac{1}{n_0}\sum_{k_1=1}^{n_0} \sigma\Big(\frac{1}{\sqrt{n_0}}\sum_{k_0=1}^{n_0} W_0(k_1,k_0)X(k_0,j_m)\Big)^2, \quad m = 1,2, \\
   \label{eq:second} \E G_1^{(n)}G_2^{(n)} & =\ind{i_1=i_2} \sum_{k_1=1}^{n_0} \frac{1}{n_0}\prod_{m=1}^2 \sigma\Big(\frac{1}{\sqrt{n_0}}\sum_{k_0=1}^{n_0} W_0(k_1,k_0)X(k_0,j_m)\Big).
 \end{align}

It therefore suffices to prove that the right-hand sides above converge in probability to $\E \sigma^2(g)$ and $\ind{i_1 =i_2} (\E \sigma(g))^2$ respectively.
Conditionally on $X$, the summands in the outer sum are independent and bounded, so by a conditional application of Chebyshev's inequality, it is enough to show the convergence in probability of the right-hand sides integrated with respect to $W_0$.  Note however that after this integration, they become
\begin{displaymath}
  \E \sigma^2(\gamma_m^{(n)}), \textrm{ $m=1,2$ and } \ind{i_1=i_2}\E \sigma(\gamma_1^{(n)})\sigma(\gamma_2^{(n)}).
\end{displaymath}
where $(\gamma_1^{(n)},\gamma_2^{(n)})$ is a centered Gaussian vector with $\E (\gamma_m^{(n)})^2 = \frac{1}{n_0}\sum_{k=1}^{n_0} X(k,j_m)^2 \stackrel{n\to \infty}{\to} 1$ a.s. and
$\E \gamma_1^{(n)}\gamma_2^{(n)} = \frac{1}{n_0}\sum_{k=1}^{n_0} X(k,j_1)X(k,j_2) \stackrel{n\to \infty}{\to}\ind{j_1=j_2}$ a.s. by the Law of Large Numbers.
Thus, $(\gamma_1^{(n)},\gamma_2^{(n)})$ converges in distribution to a centered Gaussian vector $(\gamma_1^{(\infty)},\gamma_2^{(\infty)})$ in $\R^2$, with $\E (\gamma_1^{(\infty)})^2 = \E (\gamma_2^{(\infty)})^2 = 1$ and $\E \gamma_1^{(\infty)}\gamma_1^{(\infty)} = \ind{j_1=j_2}$. Recall that we assume that $(i_1,j_1)\neq (i_2,j_2)$, so $i_1=i_2$ implies $j_1\neq j_2$, which in turn gives independence of $\gamma_1^{(\infty)}$ and $\gamma_2^{(\infty)}$. Thus, indeed  the right-hand sides of \eqref{eq:first} and \eqref{eq:second} converge in probability to $\E \sigma(g)^2$  and $\ind{i_1=i_2}(\E \sigma(g))^2$ respectively, where $g$ is a standard Gaussian variable. Together with the discussion above, this proves the claimed weak convergence of $\sqrt{an}(Z_2(i_1,j_1),Z_2(i_2,j_2))$ to $(\sigma(G_1),\sigma(G_2))$.

Having established the convergence of $\sqrt{na}(Z_2(i_1,j_1),Z_2(i_2,j_2))$ we can analyse the variance of $F(Z_2)$. To simplify the notation, we will write $A_n\simeq B_n$ if $A_n/B_n \to 1$ as $n\to \infty$. We have
\begin{displaymath}
  \E F(Z_2) = \frac{1}{n \sqrt{a}}\sum_{j=1}^n \E \sqrt{an} Z_2(1,j) = \frac{1}{\sqrt{a}}\E \sqrt{an} Z_2(1,1) \simeq \frac{1}{\sqrt{a}}\E\sigma(G_1)
\end{displaymath}
and
\begin{multline*}
  \E F(Z_2)^2 = \frac{1}{an^2 }\sum_{j=1}^n (\sqrt{an}\E Z_2(1,j))^2 + \frac{1}{a n^2 }\sum_{1\le j_1\neq j_2 \le n} \E (\sqrt{an}Z_2(1,j_1))(\sqrt{an}Z_2(1,j_2)) \\
  \simeq \frac{1}{a}\E\sigma(G_1)\sigma(G_2),
\end{multline*}
where $\E G_1^2 = \E\sigma(g)^2$ and $\Cov(G_1,G_2) = (\E \sigma(g))^2$.
Thus $\Var(F(Z_2)) \simeq \frac{1}{a} \Cov(\sigma(G_1),\sigma(G_2))$.
One can prove that $\Cov(\sigma(G_1),\sigma(G_2)) > 0$, by exploring the fact that positively correlated Gaussians are positively associated \cite{MR665603}, more precisely one can follow the proof in \cite{MR665603,krishnapur2020gaussian} and observe that the inequalities used therein in the case of the function $\sigma$ become strict on a set of positive measure. Alternatively, one may decompose $(G_1,G_2) = (\widetilde{G}_1 + \gamma,\widetilde{G_2}+\gamma)$ where $\widetilde{G_1},\widetilde{G_2},\gamma$ are independent nondegenerate Gaussian variables, $\widetilde{G}_1,\widetilde{G_2}$ identically distributed, and observe that $\Cov(\sigma(G_1),\sigma(G_2)) = \Var(f(\gamma))$, where $f(x) = \E \sigma(\widetilde{G}_1 + x)$ is strictly increasing. This ends the proof of \eqref{eq:variance-behaviour}.

\subsubsection*{Bottlenecks}
Let us conclude this discussion with an observation that one cannot recover better concentration by adding an additional independent wide layer to the network if an intermediate layer is very small. We will demonstrate it by extending the previous example, however only in the degenerate case, when $a = 1/n$, i.e., $n_2 = 1$, and only for unbounded activation function. In this case $\Var(F(Z_2))$ grows like $cn$. Let us add an additional layer of width $n_3\to \infty$, with $W_2$ being an independent $n_3\times 1$ matrix of i.i.d. standard Gaussian variables and $\sigma_2(x) = |x|$. We thus have for $i,j = 1,\ldots,n$,
\begin{displaymath}
  Z_3(i,j) = \frac{1}{\sqrt{n_3}} |W_2(i,1) Z_2(1,j)|.
\end{displaymath}

Consider now $H\colon \R^{n\times n} \to \R$, given by $H(z) = \frac{1}{\sqrt{n n_3}} \sum_{i=1}^{n_3}\sum_{j=1}^n z_{ij}$, which is $1$-Lipschitz. Using \eqref{eq:convergence-in-law} (recall that when proving \eqref{eq:convergence-in-law} we allowed $a$ to be arbitrary, also varying with $n$), we get
\begin{displaymath}
  \E H(Z_3) = \frac{1}{n^{1/2}n_3} \sum_{i=1}^{n_3}\sum_{j=1}^n \E |W_2(i,1) Z_2(1,j)| = \sqrt{n} \E |W_2(i,1)| \cdot \E|Z_2(1,1)| \simeq \sqrt{\frac{2n}{\pi}} \E \sigma(G_1)
\end{displaymath}
and
\begin{multline*}
  \E H(Z_3)^2 = \frac{1}{n n_3^2}\sum_{i=1}^{n_3}\sum_{j=1}^n |W_2(i,1)|^2 |Z_2(1,j)|^2 \\
  + \frac{1}{n n_3^2}\sum_{\stackrel{(i_1,j_1)\neq (i_2,j_2)}{1\le i_1,i_2\le n_3, 1\le j_1,j_2 \le n}} \E |W_2(i_1,1) W_2(i_2,1)| |Z_2(1,j_1)Z_2(1,j_2)|\\
  \simeq \frac{n_3 n(n-1)}{n n_3^2}\E Z_2(1,1)\E Z_2(1,2) + \frac{n n_3(n_3-1)}{n n_3^2} \frac{2}{\pi} \E Z_2(1,1)^2 \\
  + \frac{n(n-1)n_3(n_3-1)}{n n_3^2} \frac{2}{\pi}\E Z_2(1,1)Z_2(1,2)\\
\simeq \frac{2}{\pi} n \E \sigma(G_1)\sigma(G_2).
\end{multline*}
Thus $\Var H(Z_3) \simeq \frac{2}{\pi}n \Cov(\sigma(G_1),\sigma(G_2))$, i.e., it again grows linearly with $n$, which means that $Z_3$ does not satisfy dimension free concentration.

\subsubsection{Dependence on $\eta$} Let us now show by a simple example that if the data dimension $n_0$ is small, then even if $\alpha$ and $\beta$ tend to zero, the number of samples $n$ tends to $\infty$, and there is no bias, the subgaussian concentration property of $Z_L$ does not improve, i.e., the constant $c$ in \eqref{eq:bounded-activation-concentration-simplified-subgaussian} remains bounded from above. Observe that the example considered in the previous subsection, and in particular the inequality \eqref{eq:variance-behaviour}, shows an analogous phenomenon for $\alpha,\beta$ fixed and $n,n_0\to \infty$. Thus, to get arbitrary large constants $c$ in the concentration inequalities of Theorems \ref{thm:bounded-activation}, \ref{thm:unbounded-activation}, and \ref{thm:unbounded-activation-independence} we indeed need to have $\alpha\to 0, \eta\to \infty$ simultaneously.

Let us simply consider $L = 1$ with $B_0 = 0$, $n_0 = 1$, $X = (X(1),\ldots,X(n))$, $W_0 = (W_0(1),\ldots,W_0(n_1))^T$ where $X(j),W_0(i)$'s are i.i.d. random variables, uniform on $(-1,1)$ and let  $\sigma_0(x) = \min(1,|x|)$. Consider $F\colon R^{n_1\times n} \to \R$ given by $F(z) = \frac{1}{\sqrt{n_1n}} \sum_{i=1}^{n_1}\sum_{j=1}^n|z_{ij}|$, which is $1$-Lipschitz with respect to the Hilbert--Schmidt norm. Note that $\sigma_0(W_0(i)X(j)) = |W_0(i)X(j)|$.

Taking into account independence and the equality $\E |X(j)| = \E|W_0(i)| = 1/2$, we get  $\E F(Z_1) = \frac{1}{n_1\sqrt{n}} \sum_{i=1}^{n_1}\sum_{j=1}^n \E |W_i|\E|X_j| = \sqrt{n}/4$.

By the CLT we have $\lim_{n_1\to \infty} \p(\sum_{i=1}^{n_1} |W_0(i)| \ge n_1/2) = 1/2$, moreover for any fixed $t > 1$,
\begin{displaymath}
   \liminf_{n\to \infty} \p\Big(n^{-1/2}\sum_{j=1}^{n} |X(j)| \ge \frac{\sqrt{n}}{2} + t\Big) \ge \exp(-C t^2)
\end{displaymath}
for some absolute constant $C$. By independence, this implies that for large $n_1,n$, irrespectively of the  ratio $n/n_1$, with probability at least $4^{-1}\exp(-2C t^2)$, we have
\begin{displaymath}
  F(Z_1) = \frac{1}{n_1}\sum_{i=1}^{n_1} |W_0(i)| \cdot \frac{1}{n^{1/2}} \sum_{j=1}^n |X(j)| \ge \frac{1}{2}\Big(\frac{\sqrt{n}}{2} + t\Big) = \E F(Z_1) + t/2.
\end{displaymath}
This shows that the best constant $c$ in \eqref{eq:bounded-activation-concentration-simplified-subgaussian} remains in this case bounded from above, even if one takes $n_1$ so large that Assumption A1 is satisfied with $\alpha,\beta \to 0$. We remark that in this example one could alternatively estimate the variance, similarly as in the previous ones.

\subsubsection{ Dependence on $\beta$ in Theorem \ref{thm:unbounded-activation}}\phantom{a}

We will now show that in the case of dependent matrices, the constant in inequality \eqref{eq:unbounded-activation-concentration} of Theorem \ref{thm:unbounded-activation} must depend on $\beta$. To this end, consider $L = 2$, $n_0, n$ arbitrary and $n_1 = n_2\beta$. Assume further that $B_0 = B_1 = 0$, the entries of $X$ and $W_0$ are i.i.d. $N(0,1)$ random variables and $W_1 = W_0^T$. Since $(X,W_0)$ has dimension-free subgaussian concentration property and the embedding $W \mapsto (W_0,W_0^T)$ is $\sqrt{2}$-Lipschitz, $\XX$ also satisfies \eqref{eq:concentration-assumption} with an absolute constant $\rho$ (see Proposition \ref{prop:concentration-repetition}). If we set $\sigma_0(x) = \sigma_1(x) = x$, then $Z_2 = \frac{1}{\sqrt{\beta} n_2\sqrt{n_0}} W_0^TW_0 X$.

Let us denote for simplicity the columns of $W_0$ by $W_0(1),\ldots,W_0(n_0)$ (thus $W_0(i)$ are i.i.d. standard Gaussian vectors in $\R^{n_2\beta}$. Then
\begin{displaymath}
Z_2(1,1) = \frac{1}{\sqrt{\beta}n_2\sqrt{n_0}} \sum_{i=1}^{n_0} \langle W_0(1),W_0(i)\rangle X(i,1).
\end{displaymath}
By independence $\E Z_2(1,1) = 0$ and for $\beta > 2$,
\begin{multline*}
\E Z_2(1,1)^2 = \frac{1}{\beta n_2^2 n_0} \sum_{i=1}^{n^2} \E \langle W_0(1),W_0(i) \rangle^2 \ge \frac{1}{\beta n_2^2 n_0} \E \langle W_0(1),W_0(1)\rangle^2\\
=   \frac{1}{\beta n_2^2 n_0} (3 \beta n_2 + (\beta n_2)(\beta n_2 - 1)) \ge  \frac{\beta}{2 n_0}.
\end{multline*}

Thus if $\beta$ is of order greater than $n_0$, $\Var(Z_2(1,1))$ diverges, which shows that the constant $c$ in \eqref{eq:unbounded-activation-concentration} cannot be taken independent of $\beta$. Note that the remaining constants in Assumption A2, i.e., $\rho,\lambda_\sigma, \tau,\alpha, \eta$ are in this example bounded independently of $n$ if we take $n_1$ large enough with respect to $n$.

\section{Proofs of general concentration inequalities from Sections \ref{sec:bounded-activation} to \ref{sec:polynomials}}\label{sec:proofs}

In this section we will present the proofs of our main results. We will start by introducing several basic lemmas, related to the general theory of concentration of measure, which will be used in the main arguments. We remark that the core of our proofs in situations when no independence is assumed is inspired by the approach developed in order to study deviation inequalities for polynomials in independent random variables with coefficients in a Banach space. This approach was used by Borell  \cite{borell1984taylor}, Arcones and Gin\'e \cite{MR1201060}, and Ledoux \cite{MR1071528} (see also \cite{MR1102015}). Originally they were focused on the deviation above some multiple of the median, rather than strict concentration around the median or mean, however one can modify them to yield a stronger statement. Recently, similar arguments were also employed by Louart outside the polynomial case to give concentration for general smooth functions with bounded derivatives of higher order \cite{louart2025operationconcentrationinequalities}. In general such inequalities yield concentration estimates in terms of certain expected or maximal operator norms of higher order derivatives. In our setting we assume just the Lipschitz property of the activation functions, so we need to modify the original approach to this setting. This is possible thanks to the structure of the neural network and the possibility of bounding certain operator norms of random matrices along the induction scheme.

The proofs under additional independence assumptions also rely on induction and involve bounding the operator norms, however formally they resemble more the approach used, e.g.,  in \cite{MR3383337,MR4181362} to obtain concentration inequalities expressed in terms of certain injective tensor product norms of higher order derivatives (different than in the case of inequalities considered by Louart). Again, the special structure of the neural network allows to perform additional analysis during the induction, which leads  to bounds valid just under Lipschitz continuity.

\subsection{Auxiliary lemmas}

For $A \subseteq \R^N$, set
\begin{displaymath}
  A_t  = \{x\colon \dist(x,A) < t\}
\end{displaymath}

We will use the following lemma which is quite standard and a part of folklore (a similar statement in a slightly different formalism can be found, e.g., in \cite[Lemma 1.1]{MR1849347}). For completeness we provide its proof (as well as proofs of other auxiliary lemmas) in the Appendix.

\begin{lemma}\label{le:enlargement}
If a random vector $\XX$ satisfies \eqref{eq:concentration-assumption}, then for every $p \in (0,1/2]$, every set $A \subseteq \R^N$ with $\p(\XX \in A) \ge p$,
and every $t > 0$,
\begin{displaymath}
  \p(X \not\in A_t) \le \frac{2}{p}\exp\Big(-\frac{t^2}{4\rho^2}\Big).
\end{displaymath}
\end{lemma}

The next, well known lemma allows to pass from concentration around median to concentration around mean.
\begin{lemma}\label{le:concentration-around-mean}
Assume that a random vector $Y$ with values in $\R^N$ satisfies the subgaussian concentration property \eqref{eq:concentration-assumption}. Then for every $1$-Lipschitz function $F\colon \R^N \to \R$ and every $t > 0$,
\begin{displaymath}
  \p(|F(Y) - \E F(Y)| \ge t) \le 2\exp\Big(-\frac{t^2}{8\rho^2}\Big).
\end{displaymath}
\end{lemma}

Let us also formulate a less precise but more general lemma concerning equivalence of concentration around the median and the mean for more general profiles, which appear in our main results.

\begin{lemma}\label{le:concentration-equivalence}
Let $Y$ be a random vector with values in $\R^N$. For arbitrary $\alpha_1,\ldots,\alpha_m > 0$ the following conditions are equivalent:

a) There exist positive constants $a_1,\ldots,a_m$, such that for every 1-Lipschitz function $F \colon \R^N \to \R$, any median $\Med F(Y)$, and any $t > 0$,
\begin{align}\label{eq:conc-median}
  \p(|F(Y) - \Med F(Y)| \ge t) \le 2\exp(-\min(a_1 t^{\alpha_1},\ldots,a_m t^{\alpha_m})),
\end{align}

b) There exist positive constants $b_1,\ldots,b_m$, such that for every 1-Lipschitz function $F \colon \R^N \to \R$,
\begin{align}\label{eq:conc-mean}
  \p(|F(Y) - \E F(Y)| \ge t) \le 2\exp(-\min(b_1 t^{\alpha_1},\ldots,b_m t^{\alpha_m})).
\end{align}
More precisely, there exists $c > 0$, depending only on the exponents $\alpha_1,\ldots,\alpha_m$, such that if a) holds, then so does b) with $b_i = c\alpha_i$, and if b) holds, then so does a) with $a_i = cb_i$.
\end{lemma}

We will also need simple bounds on the deviation and expected value of the operator norm of centered, subgaussian random matrices. They are well-known (see for instance \cite[Theorem 4.4.5]{MR3837109}) and follow from standard epsilon-net estimates, we provide them in a lemma below. For completeness, in the Appendix we provide a full proof. We do not try to optimize constants in the estimates.

\begin{lemma}\label{le:operator-norm}
Consider a $k \times m$ centered random matrix $U$, satisfying the subgaussian concentration property with constant $\rho$.
Then, for every $s \ge 0$,
\begin{displaymath}
  \p\Big(\|U\|_{op} > \frac{16}{7} s(\sqrt{k} + \sqrt{m})\Big) \le 9^{k+m}\cdot 2 \exp(-(s^2(\sqrt{k} + \sqrt{m})^2/2\rho^2).
\end{displaymath}
As a consequence
\begin{displaymath}
\E \|U\|_{op} \le 8(\sqrt{k} + \sqrt{m})\rho.
\end{displaymath}

\end{lemma}

We will also often rely on the classical Hoffman--Wielandt inequality \cite{MR52379} (see \cite[Theorem 4.2.5]{MR3113826} for the version with singular values).

\begin{theorem}\label{thm:HW} If $A, B$ are $n\times n$ Hermitian matrices and $\gamma_i(A)$, $\gamma_i(B)$, $i \le n$, the eigenvalues of $A$, $B$ respectively, arranged in non-increasing order, then
\begin{displaymath}
  \sum_{i=1}^n (\gamma_i(A) - \gamma_i(B))^2 \le \|A - B\|_{HS}^2.
\end{displaymath}
The same inequality holds if $A,B$ are $m\times n$ matrices and $\gamma_i(A), \gamma_i(B)$ are their singular values arranged in non-increasing order.
\end{theorem}

Finally, we will need the following lemma, which we will repeatedly use to pass from tails to moments of concentrated random variables.

\begin{lemma}\label{le:moments} For $r > 0, p > 0$,
\begin{displaymath}
  \int_0^\infty p t^{p-1}\exp(-c t^r)dt = \int_0^\infty  p/r (u/c)^{(p-1)/r - (r-1)/r} c^{-1} \exp(-u)du = \frac{p}{r c^{p/r}} \Gamma(p/r).
\end{displaymath}
\end{lemma}

\subsection{Bounded and Lipschitz activation function}

We will now prove Theorem \ref{thm:bounded-activation}.

\begin{proof}[Proof of Theorem \ref{thm:bounded-activation}]

In the proof the constants $C_i$ are allowed to depend on $L,\lambda_\sigma,\kappa, \beta$ (note that they do not depend on $\alpha$, $\eta$, $\rho$, the dimensions $n, n_\ell$ or $t$). If we write, e.g., $C_i(a)$, it means that they additionally depend on the parameter $a$.

For $L = 0$, the statement follows immediately from the concentration assumption (with $c = \eta/\rho^2$), so in what follows we will assume that $L \ge 1$.

Consider first deterministic data
\begin{align}\label{eq:deterministic-data}
\xx = (x,(w_\ell,b_\ell)_{\ell=0,\ldots,L-1}), \quad \hxx = (\hx,(\hw_\ell,\hb_\ell)_{\ell=0,\ldots,L-1}),
\end{align}
where $x,\hx$ are $n_0\times n$ matrices, $w_\ell, \hw_\ell$ are $n_{\ell+1}\times n_\ell$ matrices and $b_\ell,\hb_\ell \in \R^{n_{\ell+1}}$.

Define
\begin{align}\label{eq:z-l-recurrence}
z_0 & = \frac{x}{\sqrt{n_0}},\; z_{\ell+1} =   \frac{1}{\sqrt{n_{\ell+1}}}\sigma_\ell\Big(w_\ell z_\ell + b_\ell \1_{n}^T \Big),\nonumber\\
\hz_0 &= \frac{\hx}{\sqrt{n_0}},\; \hz_{\ell+1} =   \frac{1}{\sqrt{n_{\ell+1}}}\sigma\Big(\hw_\ell \hz_\ell + \hb_\ell\1_{n} ^T\Big).
\end{align}

Thus $z_\ell,\hz_\ell$ are the matrices describing the evolution of the initial data  $x,\hx$ over the network with weights $w_\ell,b_\ell$ and $\hw_\ell,\hb_\ell$ respectively.

Assume that for some number $M$ the following inequalities hold:
\begin{align}\label{eq:iteration-conditions}
\|\xx - \yy\|_{2} & \le t, \nonumber\\
\|\hw_\ell\|_{op} & \le M(\sqrt{n_{\ell+1}} + \sqrt{n_\ell}) \textrm{ for $\ell = 0,\ldots,{L-1}$}, \nonumber \\
\|\hx\|_{op} &\le M(\sqrt{n} + \sqrt{n_0}).
\end{align}

For $\ell = 1,\ldots,L-1$ we have
\begin{align*}
\|z_{\ell+1} - \hz_{\ell+1}\|_{HS} \le&  \frac{\lambda_\sigma}{\sqrt{n_{\ell+1}}} \|w_{\ell} z_{\ell}  + b_{\ell}\1_{n}^T - \hw_{\ell} \hz_{\ell}  -  \hb_{\ell}\1_{n}^T\|_{HS}\\
\le & \frac{\lambda_\sigma}{\sqrt{n_{\ell+1}}} \|w_{\ell} z_{\ell}   - \hw_{\ell} \hz_{\ell}  \|_{HS}
+ \frac{\lambda_\sigma}{\sqrt{n_{\ell+1}}} \| b_{\ell}\1_{n}^T - \hb_{\ell}\1_{n} ^T\|_{HS}\\
\le & \frac{\lambda_\sigma}{\sqrt{n_{\ell+1}}} \|w_{\ell} - \hw_{\ell}\|_{HS} \|z_{\ell}\|_{op} + \frac{\lambda_\sigma}{\sqrt{n_{\ell+1}}} \|\hw_{\ell}\|_{op}\|z_{\ell} - \hz_{\ell}\|_{HS}\\
& + \lambda_\sigma\frac{\sqrt{n}}{\sqrt{n_{\ell+1}}}\|b_{\ell} - \hb_{\ell}\|_2\\
\le & \frac{\lambda_\sigma}{\sqrt{n_{\ell+1}}} t \frac{\sqrt{n n_{\ell}}}{\sqrt{n_\ell}}\kappa + M\frac{\lambda_\sigma}{\sqrt{n_{\ell+1}}} (\sqrt{n_{\ell+1}} + \sqrt{n_{\ell}}) \|z_{\ell} - \hz_{\ell}\|_{HS}\\
&+ \lambda_\sigma\frac{\sqrt{n}}{\sqrt{n_{\ell+1}}}t\\
\le & M\lambda_\sigma(1 + \sqrt{\beta})\|z_{\ell} - \hz_{\ell}\|_{HS} + \lambda_\sigma \sqrt{\alpha} (1+ \kappa)t,
\end{align*}
where in the fourth inequality we used the trivial bound $\|z_\ell\|_{op} \le \|z_\ell\|_{HS} \le \sqrt{nn_\ell}\frac{\kappa}{\sqrt{n_{\ell}}}$.

By induction,
\begin{align}\label{eq:bound-by-z_1}
  \|z_L - \hz_L\|_{HS} \le &  C_1(M)\sqrt{\alpha}t  + C_2(M)\|z_1 - \hz_1\|_{HS}
\end{align}
(recall that $C_i(M)$ depends only on $L,\lambda_\sigma, \kappa,\beta$ and $M$).

To bound $\|z_1 - \hz_1\|_{HS}$ we proceed similarly.We have, however, to handle this term separately, since we cannot use the $L_\infty$ bounds on $\sigma_\ell$ to control $\|z_0\|_{op}$. Instead we will rely on the last inequality of \eqref{eq:iteration-conditions}. We have

\begin{align*}
  \|z_{1} - \hz_{1}\|_{HS} \le&  \frac{\lambda_\sigma}{\sqrt{n_{1}}} \|w_{0} z_{0}  +  b_{0}\1_{n}^T - \hw_{0} \hz_{0}  - \hb_{0}\1_{n}^T\|_{HS}\\
\le & \frac{\lambda_\sigma}{\sqrt{n_{1}}} \|w_{0} z_{0}   - \hw_{0} \hz_{0}  \|_{HS}
+ \frac{\lambda_\sigma}{\sqrt{n_{1}}} \|  b_{0}\1_{n}^T - \hb_{0}\1_{n} ^T\|_{HS}\\
\le& \frac{\lambda_\sigma}{\sqrt{n_{1}}} \|(w_0 - \hw_0)(z_0 - \hz_0)\|_{HS} + \frac{\lambda_\sigma}{\sqrt{n_{1}}} \|\hw_0(z_0-\hz_0)\|_{HS} \\
& + \frac{\lambda_\sigma}{\sqrt{n_{1}}} \|(w_0 - \hw_0)\hz_0\|_{HS}
+ \frac{\lambda_\sigma\sqrt{n}}{\sqrt{n_{1}}} \|b_0 - \hb_0\|_2\\
\le & \frac{\lambda_\sigma}{\sqrt{n_{1}}} \|w_0-\hw_0\|_{HS}\|z_0 - \hz_0\|_{HS} +  \frac{\lambda_\sigma}{\sqrt{n_{1}}} \|\hw_0\|_{op}\|z_0-\hz_0\|_{HS}\\
& + \frac{\lambda_\sigma}{\sqrt{n_{1}}} \|w_0 - \hw_0\|_{HS}\|\hz_0\|_{op} + \frac{\lambda_\sigma\sqrt{n}}{\sqrt{n_{1}}} \|b_0 - \hb_0\|_2\\
\le & \frac{\lambda_\sigma}{\sqrt{n_{1}}} \frac{t^2}{\sqrt{n_0}} + \frac{\lambda_\sigma}{\sqrt{n_{1}}} M(\sqrt{n_0} + \sqrt{n_1})\frac{t}{\sqrt{n_0}}\\
&+ \frac{\lambda_\sigma}{\sqrt{n_{1}}} t\frac{M(\sqrt{n_0} + \sqrt{n})}{\sqrt{n_0}} + \frac{\lambda_\sigma\sqrt{n}}{\sqrt{n_{1}}} t\\
\le & \frac{\lambda_\sigma}{\sqrt{nn_0}}\sqrt{\alpha} t^2 + \frac{\lambda_\sigma}{\sqrt{n_0}} M(\sqrt{\beta}+1)t + \frac{\lambda_\sigma}{\sqrt{n_0}}M(\sqrt{\beta}+\sqrt{\alpha}) t + \lambda_\sigma \sqrt{\alpha} t.
\end{align*}

Combining this with the inequality \eqref{eq:bound-by-z_1} (note the remark following it) and taking into account the Lipschitz property of $F$ together with the assumption $n_0 \ge \eta$, we obtain
\begin{align}\label{eq:deterministic-conclusion}
|F(z_L) - F(\hz_L)| \le  \|z_L - \hz_L\|_{HS} \le \Big(C_3(M,\alpha,\eta) t + \frac{C_4(M,\alpha,\eta)}{\sqrt{n}} t^2\Big).
\end{align}
where for $i=3,4$ and the parameters $L,\lambda_\sigma,\kappa, \beta,M$ fixed,
\begin{align}\label{eq:it-goes-to-zero}
\lim_{\alpha \to 0, \eta \to \infty} C_i(M,\alpha,\eta) = 0.
\end{align}
Let $M$ be such that
\begin{align}\label{eq:condition-on-M}
\max_{0 \le \ell \le L-1}\p(\|W_\ell\|_{op} \ge M(\sqrt{n_\ell} + \sqrt{n_{\ell+1}})), \p(\|X\|_{op} \ge M (\sqrt{n} + \sqrt{n_0})) \le \frac{1}{4(L+1)}
\end{align}
and let $\mathcal{A}\subset \R^N$ be the set of all $\hxx = (\hx,(\hw_\ell,\hb_\ell)_{\ell=0,\ldots,L-1})$ such that
\begin{align*}
  F(\hz_L) &\le \Med F(Z_L),\\
  \|\hw_\ell\|_{op} & \le M(\sqrt{n_\ell} + \sqrt{n_{\ell+1}}), \; \ell = 0,\ldots,L-1,\\
  \|\hx\|_{op} & \le M(\sqrt{n} + \sqrt{n_0}).
\end{align*}

By the definition of the median, and the union bound, we have $\p(\XX \in \mathcal{A}) \ge 1/4$. Thus, by Lemma \ref{le:enlargement}, for all $t \ge 0$,
\begin{displaymath}
\p(\XX \notin \mathcal{A}_{t}) \le 8 \exp(-t^2/4\rho).
\end{displaymath}

Note that if $\XX \in \mathcal{A}_t$, then we can find $\hxx = (\hx,(\hw_\ell,\hb_\ell)_{\ell=0,\ldots,L-1})$ such that if we set $\xx = \XX$, the conditions \eqref{eq:z-l-recurrence} are satisfied and moreover $F(\hz_L) \le \Med F(Z_L)$. In particular, the inequality \eqref{eq:deterministic-conclusion} holds with $Z_L$, instead of $z_L$.

As a consequence with probability at least $1 - 8\exp(-t^2/4\rho^2)$, we have
\begin{displaymath}
  F(Z_L) - \Med F(Z_L) \le F(Z_L) - F(\hz_L) \le \Big(C_3(M,\alpha,\eta) t + \frac{C_4(M,\alpha,\eta)}{\sqrt{n}} t^2\Big).
\end{displaymath}
Applying the same argument to $-F$ and using the union bound, we get
\begin{displaymath}
  \p\Big(|F(Z_L) - \Med F(Z_L)| \ge \Big(C_3(M,\alpha,\eta) t + \frac{C_4(M,\alpha,\eta)}{\sqrt{n}} t^2\Big)\Big) \le 16\exp(-t^2/4\rho^2),
\end{displaymath}
which by a substitution $t := \min(t/2 C_3(M,\alpha,\eta),  n^{1/4}(t/2C_4(M,\alpha,\eta))^{1/2})$, implies that
\begin{align}\label{eq:in-terms-of-M}
  \p(|F(Z_L) - \Med F(Z_L)| \ge t) \le 16\exp\Big(-\frac{1}{4\rho^2}\min\Big(\frac{t^2}{4C_3(M,\alpha,\eta)^2 },\frac{\sqrt{n} t}{2 C_4(M,\alpha,\eta)}\Big)\Big).
\end{align}

We will now show that one can find $M$, satisfying \eqref{eq:condition-on-M} and depending only on $\rho$ and $L$.

By Lemma \ref{le:operator-norm}, for any centered $k \times m$ random matrix satisfying the subgaussian concentration inequality with parameter $\rho$,
\begin{displaymath}
\p(\|U\| > M(\sqrt{k} + \sqrt{m})) \le \frac{\E \|U\|}{M(\sqrt{k}+\sqrt{m})} \le \frac{8\rho}{M}.
\end{displaymath}
Thus it is indeed enough to take $M = 32\rho (L+1)$ (see also Remark \ref{re:constant-net} below).

From \eqref{eq:in-terms-of-M} and \eqref{eq:it-goes-to-zero} it follows that there exists a constant $c$, depending on $L,\rho,\lambda_\sigma,\kappa,\beta,\alpha,\eta$, such that $c \to \infty$, when $\alpha \to 0$ and $\eta \to \infty$ simultaneously, while the other parameters are fixed, and moreover whenever the Assumption A1 is satisfied, then for any $t \ge 0$,
\begin{align}\label{eq:two-level}
  \p(|F(Z_L) - \Med F(Z_L)| \ge t) \le 16\exp\bigl(-c \min(t^2,\sqrt{n} t)\bigr).
\end{align}

To prove \eqref{eq:bounded-activation-concentration-simplified-subgaussian} it is enough to observe that in the range when the left-hand side of the above inequality is nontrivial, the second term in the exponent actually dominates (up to a constant) the first one. Indeed, consider any $\xx = (x,(w_\ell,b_\ell)_{\ell=0,\ldots,L-1})$ in the support of $\XX$ and let $z_L$ be defined through the recurrence \eqref{eq:z-l-recurrence}. Then
\begin{displaymath}
  |F(Z_L) - F(z_L)| \le \|Z_L - z_L\|_{HS} \le \frac{1}{\sqrt{n_L}} 2\kappa \sqrt{n_L n} = 2\kappa \sqrt{n}.
\end{displaymath}
Thus trivially for $t > 2\kappa\sqrt{n}$, $ \p(|F(Z_L) - \Med F(Z_L)| \ge t) = 0$. On the other hand, for $t \le 2\kappa \sqrt{n}$, we have
\begin{displaymath}
\sqrt{n} t \ge \frac{t^2}{2\kappa}
\end{displaymath}
and thus \eqref{eq:two-level} gives
\begin{displaymath}
  \p(|F(Z_L) - \Med F(Z_L)| \ge t) \le 16\exp\bigl(-c \min(t^2,t^2/(2\kappa)\bigr).
\end{displaymath}
To finish the proof it is now enough to adjust the constants.
\end{proof}

\begin{remark}\label{re:constant-net} Since we do not provide explicit dependence of constants on the parameters, when choosing the constant $M$ in the above argument, we used just Markov's inequality in combination with the second estimate of Lemma \ref{le:operator-norm}. In fact, the proof of Theorem \ref{thm:bounded-activation} leads to explicit constants, which however are given by a rather complicated expression involving sums of some geometric series. If one is interested in obtaining an explicit estimate on the constants, it is better to choose $M$ with help of the first estimate of Lemma \ref{le:operator-norm}. The dependence of $M$ on $L$ is then logarithmic, instead of the linear one obtained above (still the remaining parts of the proof give much worse dependence of the constants on $L$, which is the main reason why we do not make this dependence explicit).
\end{remark}

\subsection{General Lipschitz activation function}

In this section we will prove Theorems \ref{thm:unbounded-activation} and \ref{thm:unbounded-activation-independence}.

\begin{proof}[Proof of Theorem \ref{thm:unbounded-activation} and Proposition \ref{prop:unbounded-activation}]

Let us start with the proof of \eqref{eq:unbounded-activation-concentration}.
We will prove by induction over $L\ge 0$, that there exist constants $c_L, C_L$, depending only on $L, \rho, \lambda_\sigma, \tau, \alpha,\beta,\eta$ (in particular independent of $n$), such that

\medskip

\noindent \textbullet\,

\begin{align}\label{eq:it-goes-tozero-induction}
\textrm{for fixed $L, \rho, \lambda_\sigma, \tau,\beta$, } \lim_{\alpha\to 0,\eta\to \infty} c_L =  \infty,
\end{align}

\noindent \textbullet\, for every random vector $\XX$, satisfying the Assumption A2,
\begin{align}\label{eq:induction-assumption-HS-norm}
\E \|Z_L\|_{HS} \le C_L \sqrt{n}
\end{align}
and for all 1-Lipschitz functions $F \colon \R^{n_\ell\times n} \to \R$, and all $t > 0$,
\begin{align}\label{eq:induction-assumption-tail}
\p(|F(Z_L) - \Med F(Z_L)| \ge t) & \le 2\exp\Big(-c_L\min_{1\le k \le L+1} \Big(n^{(k-1)/k}t^{2/k}\Big) \Big).
\end{align}

In the proof we will use two types of auxiliary constants. The constants of the form $C,C'$, etc. are allowed to depend only on $L, \rho, \lambda_\sigma, \tau, \alpha,\beta,\eta$. On the other hand, the constants denoted with the letter $D$, e.g., $D(a_1,a_2,\dots,a_m), D'(a_1,\ldots,a_m)$, etc., are allowed to depend only on the parameters $a_1,\ldots,a_m$, listed as their arguments.

\medskip

\paragraph{\bf 1. The case $\mathbf{L = 0}$} We have $\E X_{ij} = 0$ and $X_{ij}$ satisfies the subgaussian concentration inequality with parameter $\rho$, so by Lemma \ref{le:concentration-around-mean},
\begin{displaymath}
  \E X_{ij}^2 = 2 \int_0^\infty t \p(|X_{ij}| \ge t)dt \le 2\int_0^\infty t e^{-t^2/(8\rho^2)}dt = 8\rho^2.
\end{displaymath}

Thus,
\begin{displaymath}
  (\E \|Z_0\|_{HS})^2 \le \E \|Z_{0}\|_{HS}^2 = \frac{1}{n_0} \sum_{i=1}^{n_0}\sum_{j=1}^n \E|X_{ij}|^2 \le 8\rho^2 n,
\end{displaymath}
i.e., \eqref{eq:induction-assumption-HS-norm} holds with $C_0 = 2\sqrt{2}\rho$.

The inequality \eqref{eq:induction-assumption-tail} for $L = 0$ with $c_0 = \eta/\rho^2$ is immediately implied by the concentration assumption on $X$ (note that for $L = 0$ the right-hand side of \eqref{eq:induction-assumption-tail} is just $2e^{-c_0 t^2})$. The constant $c_0$ does not depend on $\alpha$ and clearly tends to $\infty$ with $\eta\to \infty$.

\paragraph*{\bf 2. The induction step} Assume that \eqref{eq:it-goes-tozero-induction}, \eqref{eq:induction-assumption-HS-norm} and \eqref{eq:induction-assumption-tail} hold for all numbers smaller than $L+1 \ge 1$. We will prove that they hold also for $L+1$.

Let us start with an estimation of the Hilbert--Schmidt norm of $Z_{L+1}$. By Lemma \ref{le:concentration-around-mean} and integration by parts we have just as for $X_{ij}$ above,
\begin{displaymath}
  \E B_L(i)^2 \le 8 \rho^2,
\end{displaymath}
so
\begin{displaymath}
\E \|B_{L}\1_n^T\|_{HS}^2 \le 8\rho^2 n n_{L+1}.
\end{displaymath}

Moreover, by Lemma \ref{le:operator-norm}, $\E \|W_L\|_{op} \le 8(\sqrt{n_{L}} + \sqrt{n_{L+1}})\rho$. Since the operator norm is $1$-Lipschitz, we have
\begin{multline*}
\E \|W_L\|_{op}^2 \le 2 (\E \|W_L\|_{op})^2 + 2\E(\|W_L\|_{op} - \E \|W_L\|_{op})_+^2 \\
\le 128(\sqrt{n_{L}} + \sqrt{n_{L+1}})^2\rho^2 + 4 \int_0^\infty t \p(\|W_L\|_{op} - \E \|W_L\|_{op} \ge t)dt\\
\le 128(\sqrt{n_{L}} + \sqrt{n_{L+1}})^2\rho^2  + 8 \int_0^\infty t \exp(-t^2/8\rho^2)dt\\
= 128(\sqrt{n_{L}} + \sqrt{n_{L+1}})^2\rho^2 + 32\rho^2 \le 160(\sqrt{n_{L}} + \sqrt{n_{L+1}})^2\rho^2.
\end{multline*}

Similarly, using the induction assumptions \eqref{eq:induction-assumption-HS-norm} and \eqref{eq:induction-assumption-tail}, we obtain
\begin{multline*}
\E \|Z_L\|_{HS}^2 \le 2 (2\E \|Z_L\|_{HS})^2 + 2\E(\|Z_L\|_{HS} - 2\E \|Z_L\|_{HS})_+^2 \\
\le 8C_L^2 n + 8\int_0^\infty t \exp(-c_L \min(t^2,t^{2/(L+1)}))dt \le C' n,
\end{multline*}
where we used that $2\E \|Z_L\|_{HS} \ge \Med \|Z_L\|_{HS}$. Note that when estimating the tail probability we ignored the dependence on $n$ of the expressions corresponding to $k >  1$ on the right hand side of \eqref{eq:induction-assumption-tail}.

As a consequence of the last three estimates, we get
\begin{multline*}
\E\|W_L Z_L + \1_nB_L^T\|_{HS} \le \E\|W_L\|_{op} \|Z_L\|_{HS} + \E \|\1_nB_L^T\|_{HS} \\
\le (\E \|W_L\|_{op}^2)^{1/2}(\E \|Z_L\|_{HS}^2)^{1/2} + \E \|\1_nB_L^T\|_{HS} \\
\le C''\bigl((\sqrt{n_L} + \sqrt{n_{L+1}}) \sqrt{n} + \sqrt{n n_{L+1}}\bigr)
\end{multline*}
and so
\begin{multline*}
  \E\| Z_L\|_{HS} \le \tau\sqrt{n} + \frac{\lambda_\sigma}{\sqrt{n_{L+1}}} \E\|W_L Z_L + \1_nB_L^T\|_{HS}\\
  \le  \tau\sqrt{n} + C''\frac{\lambda_\sigma}{\sqrt{n_{L+1}}} \bigl((\sqrt{n_L} + \sqrt{n_{L+1}}) \sqrt{n} + \sqrt{n n_{L+1}}\bigr)\\
  \le \bigl (\tau+ C'' \lambda_\sigma(\sqrt{\beta} + 2)\bigr )\sqrt{n} =: C_{L+1}\sqrt{n}.
\end{multline*}

This establishes \eqref{eq:induction-assumption-HS-norm} for $L+1$. We may now pass to the proof of the concentration inequality \eqref{eq:induction-assumption-tail} for $L+1$.

We will use the notation introduced in \eqref{eq:deterministic-data} and \eqref{eq:z-l-recurrence}, again starting with an estimate for deterministic data, under assumptions which will be later proven to hold with high probability. Assume thus that for some $M,R > 0$, the vectors $\xx = (x,(w_\ell,b_\ell)_{\ell=0,\ldots,L}), \, \hxx = (\hx,(\hw_\ell,\hb_\ell)_{\ell=0,\ldots,L})$ satisfy

\begin{align}\label{eq:iteration-conditions-general-function}
\|\xx - \hxx\|_{2} & \le t, \nonumber\\
\|\hw_\ell\|_{op} & \le M(\sqrt{n_{\ell+1}} + \sqrt{n_\ell}) \textrm{ for $\ell= 0,\ldots,{L}$}, \nonumber \\
\|\hz_\ell\|_{op} &\le R\sqrt{n} \textrm{ for $\ell = 1,\ldots,{L}$}, \\
\|\hx\|_{op} & \le R(\sqrt{n_0} + \sqrt{n}). \nonumber
\end{align}

We will use a similar decomposition as in the proof of Theorem \ref{thm:bounded-activation}, we will however have to introduce some changes, since we do not control the norm of $z_i$ (the entries of the matrix are not bounded) but the norm of $\hz_i$ (thanks to \eqref{eq:iteration-conditions-general-function}). We have for $1\le \ell \le L$,

\begin{align*}
\|z_{\ell+1} - \hz_{\ell+1}\|_{HS} \le&  \frac{\lambda_\sigma}{\sqrt{n_{\ell+1}}} \|w_{\ell} z_{\ell}  +  b_{\ell}\1_{n}^T - \hw_{\ell} \hz_{\ell}  - \hb_{\ell}\1_{n} ^T\|_{HS}\\
\le & \frac{\lambda_\sigma}{\sqrt{n_{\ell+1}}} \|w_{\ell} z_{\ell}   - \hw_{\ell} \hz_{\ell}  \|_{HS}
+ \frac{\lambda_\sigma}{\sqrt{n_{\ell+1}}} \| b_{\ell}\1_{n} ^T - \hb_{\ell}\1_{n}^T\|_{HS}\\
\le & \frac{\lambda_\sigma}{\sqrt{n_{\ell+1}}} \Big(\|w_{\ell} - \hw_{\ell}\|_{HS} \|z_{\ell}-\hz_\ell\|_{HS} + \|w_\ell-\hw_\ell\|_{HS}\|\hz_\ell\|_{op} + \|\hw_{\ell}\|_{op}\|z_{\ell} - \hz_{\ell}\|_{HS}\Big)\\
& + \lambda_\sigma\frac{\sqrt{n}}{\sqrt{n_{\ell+1}}}\|b_{\ell} - \hb_{\ell}\|_2\\
\le & \frac{\lambda_\sigma}{\sqrt{n_{\ell+1}}} t\|z_\ell-\hz_\ell\|_{HS} + \frac{\lambda_\sigma}{\sqrt{n_{\ell+1}}}R\sqrt{n}t + M\frac{\lambda_\sigma}{\sqrt{n_{\ell+1}}} (\sqrt{n_{\ell+1}} + \sqrt{n_{\ell}}) \|z_\ell - \hz_{\ell}\|_{HS}\\
&+ \lambda_\sigma\frac{\sqrt{n}}{\sqrt{n_{\ell+1}}}t\\
\le & \lambda_\sigma \Big(M(1 + \sqrt{\beta})+\frac{\sqrt{\alpha}}{\sqrt{n}}t\Big)\|z_{\ell} - \hz_{\ell}\|_{HS} + \lambda_\sigma (R+1)\sqrt{\alpha}t.
\end{align*}
By induction over $\ell$ this gives
\begin{align}\label{eq:added-later}
  \|z_{L+1} - \hz_{L+1}\|_{HS} \le D'(L,\lambda_\sigma,\alpha,\beta,M,R) \|z_1 - \hz_1\|_{HS} \sum_{k=0}^{L} \frac{t^k}{n^{k/2}} + D''(L,\lambda_\sigma,\alpha,\beta,M,R)\sum_{k=1}^{L}\frac{t^k}{n^{(k-1)/2}},
\end{align}
where for fixed values of $L,\lambda_\sigma,\beta,M,R$, we have
\begin{align*}
\lim_{\alpha\to 0} D''(L,\lambda_\sigma,\alpha,\beta,M,R) = 0.
\end{align*}

Using the same argument as above for $\ell = 0$, but estimating $\|\hz_0\|_{op} \le n_0^{-1/2}\|x\|_{op} \le Rn_0^{-1/2}(\sqrt{n_0} + \sqrt{n})$, we get
\begin{align*}
  \|z_1 - \hz_1\|_{HS} & \le \lambda_\sigma \Big(M(1 + \sqrt{\beta})+\frac{\sqrt{\alpha}}{\sqrt{n}}t\Big)\|z_{0} - \hz_{0}\|_{HS} + \lambda_\sigma (\frac{R}{\sqrt{n_0n_1}}(\sqrt{n_0}+\sqrt{n})+\sqrt{\alpha})t\\
  & \le \lambda_\sigma \Big(M(1 + \sqrt{\beta})+\frac{\sqrt{\alpha}}{\sqrt{n}}t\Big)\frac{t}{\sqrt{n_0}}
  + \lambda_\sigma (\frac{R}{\sqrt{n_0}}(\sqrt{\beta}+\sqrt{\alpha})+\sqrt{\alpha})t\\
  &\le D'''(\lambda_\sigma, \alpha,\beta,\eta,M,R)\Big(t + \frac{t^2}{\sqrt{n}}\Big),
\end{align*}
where for fixed values of $\lambda_\sigma,\beta,M,R$,
\begin{displaymath}
  \lim_{\alpha\to 0,\eta\to \infty} D''' (\lambda_\sigma, \alpha,\beta,\eta,M,R) = 0.
\end{displaymath}

Plugging this into the previous estimate, we obtain
\begin{align}\label{eq:HS-norm-of-difference-2nd-thm}
 \|z_{L+1} - \hz_{L+1}\|_{HS} \le D(L,\lambda_\sigma,\alpha,\beta,\eta,M,R)\sum_{k=1}^{L+2} \frac{t^k}{n^{(k-1)/2}},
\end{align}
with
\begin{align}\label{the-C-goes-to-zero}
\lim_{\alpha\to 0,\eta\to \infty} D(L,\lambda_\sigma,\alpha,\beta,\eta,M,R)
\end{align}
for all the parameters, except for $\alpha,\eta$, fixed (note that $D'$ in \eqref{eq:added-later} is non-decreasing in $\alpha$).

Consider now a 1-Lipschitz function $F \colon \R^{n_{\ell+1}\times n} \to \R$ and the set $\mathcal{A} \subseteq \R^N$ of all $\hxx = (\hz,(\hw_\ell,\hb_\ell)_{\ell=0,\ldots,L-1})$ such that
\begin{align*}
F(\hz_{L+1}) &\le \Med F(Z_{L+1})\\
\|\hw_\ell\|_{op} & \le M(\sqrt{n_{\ell+1}} + \sqrt{n_\ell}) \textrm{ for $\ell = 0,\ldots,L$}, \nonumber \\
\|\hz_\ell\|_{op} &\le R\sqrt{n} \textrm{ for $\ell = 1,\ldots,L$},\nonumber \\
\|\hx\|_{op} & \le R(\sqrt{n_0}+\sqrt{n}) \nonumber,
\end{align*}
where $M$ depends only on $\rho$ and $L$, while $R$ depends on $\rho,\alpha,L$ and the constants $C_\ell$, $\ell \le L$ in such a way that $\p(\XX \in \mathcal{A}) \ge 1/4$. One can choose such parameters using the union bound, concentration for $\XX$ and the induction assumption \eqref{eq:induction-assumption-HS-norm} combined with the Chebyshev's inequality. Indeed if we choose $M = 32(2L+2)\rho$, then Chebyshev's inequality and Lemma \ref{le:operator-norm} imply that
\begin{displaymath}
  \max_{i=0,\ldots,L} \p(\|W_i\|_{op} >M(\sqrt{n_{i+1}} + \sqrt{n_i}) ), \p( \|X\|_{op} >  M(\sqrt{n} + \sqrt{n_0})) \le \frac{1}{4(2L+2)}.
\end{displaymath}

Moreover, the estimate $\E\|Z_\ell\|_{op} \le \E\|Z_\ell\|_{HS}\le C_\ell\sqrt{n}$ implies that for $R \ge 4(2L+2)\max_{1\le \ell\le L} C_\ell$, we have
\begin{displaymath}
  \p(\|Z_\ell\|_{op} > R\sqrt{n}) \le \frac{1}{4(2L+2)}, \; \ell =1,\ldots,L.
\end{displaymath}
Choosing thus $R = 4(2L+2)\max(\max_{1\le \ell \le L} C_\ell,8\rho)$, we obtain by the union bound over $1 + (L+1)+L+1 = 2L+3$ events that $\p(Z\in \mathcal{A}) \ge 1 - (1/2 + (2L + 2)/(4(2L+2))) =  1/4$.

By Lemma \ref{le:enlargement}, we get $\p(\XX \in \mathcal{A}_t) \ge 1 - 8\exp(-t^2/4\rho)$. Moreover, if $\XX \in \mathcal{A}_t$, then there exists $\hxx \in \mathcal{A}$
such that $\|\xx - \hxx\|_2 \le t$ and by \eqref{eq:HS-norm-of-difference-2nd-thm}, we get
\begin{multline*}
  F(Z_{L+1}) \le F(\hz_{L+1}) + \|Z_{L+1} - \hz_{L+1}\|_{HS} \\
  \le \Med F(Z_{L+1}) + D(L,\lambda_\sigma,\alpha,\beta,\eta,M,R)  \sum_{k=1}^{L+2} \frac{t^k}{n^{(k-1)/2}}.
\end{multline*}
Applying the same argument to $-F$ we obtain that
\begin{displaymath}
  \p\Big(|F(Z_{L+1}) - \Med F(Z_{L+1})| \ge D(L,\lambda_\sigma,\alpha,\beta,\eta,M,R)  \sum_{k=1}^{L+2} \frac{t^k}{n^{(k-1)/2}}\Big) \le 16\exp(-t^2/4\rho^2),
\end{displaymath}
which by a change of variables, an application of \eqref{the-C-goes-to-zero}, and an adjustment of constants implies \eqref{eq:induction-assumption-tail} with $L+1$ instead of $L$ and a constant $c_{L+1}$ satisfying
$\lim_{\alpha\to 0, \eta\to \infty} c_{L+1} = \infty$ for each fixed choice of the remaining parameters $L,\rho,\lambda_\sigma,\tau,\beta$. Note that in order to get the convergence of the constant $c_{L+1}$ to $\infty$ we rely on obvious monotonicity in the dependence of optimal constants $C_\ell$ on $\alpha, \eta$ as well as the optimal constants $D$ in \eqref{the-C-goes-to-zero} on $R$ (the optimal constants $C_\ell$ are non-decreasing in $\alpha$ and non-increasing in $\eta$, while the optimal constant $D$ is non-decreasing in $R$). This ends the induction step and the proof of \eqref{eq:unbounded-activation-concentration}.

The inequality \eqref{eq:unbounded-but-the-last-one} follows from \eqref{eq:unbounded-activation-concentration}. Indeed, under the assumption $\|\sigma_{L-1}\|_\infty \le \kappa$, the law of $Z_L$ is supported on the Hilbert--Schmidt ball centered at $0$, with radius $\kappa\sqrt{n}$. Thus, for any 1-Lipschitz function $F$, and $t > 2\kappa\sqrt{n}$, $\p(|F(Z) - \Med F(Z)| \ge t) = 0$. On the other hand for $t \le 2\kappa{\sqrt{n}}$,
\begin{displaymath}
   n^{(k-1)/k} t^{2/k} \ge \frac{t^{2}}{(2\kappa)^{2(k-1)/k}}, \quad k=2,\ldots,L+1.
\end{displaymath}
Substituting these inequalities into \eqref{eq:unbounded-activation-concentration} and adjusting the constants gives \eqref{eq:unbounded-but-the-last-one} and ends the proof of the theorem.
\end{proof}

\begin{remark}\label{re:L1}
Let us remark that the above proof shows that if $L = 1$ and there is no bias, i.e., $B_0 = 0$, one can obtain a better tail inequality that in the general case, since one gets an improved dependence on $n_0$. The examples presented in Section \ref{sec:discussion} show that this is not the case for $L \ge 2$, even in the absence of bias. This phenomenon is related to the fact that if one does not impose additional symmetries (for instance those in Assumption 4), then $Z_0 = \frac{1}{\sqrt{n_0}}X_0$ has operator norm of smaller order than the matrices $Z_\ell$ for $\ell \ge 1$.
\end{remark}

Let us now pass to the proof of Theorem \ref{thm:unbounded-activation-independence}.

\begin{proof}[Proof of Theorem \ref{thm:unbounded-activation-independence} and Proposition \ref{prop:unbounded-activation-independence}] Before we proceed with the argument, let us give some remarks concerning its (rather simple) structure.
Just as in the case of Theorem \ref{thm:unbounded-activation} the proof will proceed by induction, we will simultaneously prove the concentration inequality and a bound from above on the Hilbert--Schmidt norm of $Z_\ell$ (used in fact just to estimate the operator norm, which in the worst case may be of the same order). The main difference is the fact that if $Z_\ell$ is independent from $W_\ell$, instead of looking at the
Lipschitz constant of the map $z \mapsto \sigma_\ell(W_\ell z + B_\ell\1_n^T)$ for fixed value of $W_\ell, B_\ell$, which may be with high probability of the order of $\sqrt{n_{\ell+1}} + \sqrt{n_\ell}$, we may consider the Lipschitz constant of the map $z \mapsto \E \sigma_\ell(W_\ell z + B_\ell\1_n^T)$, which does not involve the parameter $n_\ell$. Heuristically, one may think that under independence it is unlikely that $Z_\ell$ and $W_\ell$ are aligned in such a way that the operator norm of the product $W_\ell Z_\ell$ is close to the product of the operator norms, whereas this is possible in the general dependent situation (see the discussion concerning the parameter $\beta$ in Section \ref{sec:discussion}).

Passing to the proof, we will demonstrate by induction on $L\ge 0$ that there exist positive constants $C_L,c_L$, depending only on $L,\rho,\lambda_\sigma,\tau,\alpha,\eta$ such that

\noindent \textbullet\,
\begin{align}\label{eq:it-goes-tozero-induction-independence}
\textrm{for fixed $L, \rho, \lambda_\sigma, \tau$, } \lim_{\alpha\to 0,\eta\to \infty} c_L =  \infty,
\end{align}

\noindent \textbullet\, for all $\XX = (X,(W_\ell,B_\ell)_{\ell=0,\ldots,L-1})$, satisfying the Assumption A3,
\begin{align}\label{eq:induction-assumption-HS-independence}
\|Z_L\|_{HS} \le C_L \sqrt{n}
\end{align}
and for any $F \colon \R^{n_{L}\times n} \to \R$, $1$-Lipschitz with respect to the Hilbert--Schmidt norm, and any $t > 0$,
\begin{align}\label{eq:induction-assumption-concentration-independence}
\p(|F(Z_L) - \E F(Z_L)| \ge t) \le 2\exp\Big(-c_L\min_{1\le k \le L+1} \Big(n^{(k-1)/k}t^{2/k}\Big) \Big).
\end{align}

In the proof, without loss of generality, we will assume that $X,(W_\ell,B_\ell)$ $\ell=0,\ldots,L$ are defined as coordinates on a product probability space. We will denote by $\E_\ell$ integration with respect to $(W_\ell,B_\ell)$ (with the other coordinates fixed) and by $\E_{\le \ell}$ integration with respect to $X,(W_i,X_i)$, $i \le \ell$. By the Fubini theorem we have $\E_\ell \E_{\le \ell - 1} = \E_{\le \ell}$. We will use a similar convention concerning conditional probabilities.

In what follows, in addition to the constants $C_L$ and $c_L$, we will use constants $D(a_1,\ldots,a_m)$, $D'(a_1,\ldots,a_m)$, etc., which are allowed to depend only on the parameters $a_1,\ldots,a_m$ (in particular if the list is empty, the corresponding constant is universal).

\paragraph{\bf 1. The case $\mathbf{L = 0}$}

This follows in almost the same way as in the proof of Theorem \ref{thm:unbounded-activation}. The bound on the operator norm follows by concentration and integration by parts, while the concentration inequality  with $c_0 = \eta/(8\rho^2)$ follows from the assumption on concentration properties of $X$ and Lemma \ref{le:concentration-around-mean}. Note that $c_0$ depends only on $\rho$ and $\eta$ and for fixed $\rho$, diverges to $\infty$ with $\eta \to \infty$.

\paragraph*{\bf 2. The induction step}

Assume that there exist constants $c_L,C_L$ satisfying \eqref{eq:it-goes-tozero-induction-independence}, \eqref{eq:induction-assumption-HS-independence} and \eqref{eq:induction-assumption-concentration-independence}. We will show that this is also the case for $L+1$.

Let $\XX = (X,(W_\ell,B_\ell)_{\ell=0,\ldots,L})$ be any vector satisfying the Assumption A3 (with $L$ replaced by $L+1$, note that $Z_{L+1}$ is a function of $\XX$).

We will start with \eqref{eq:induction-assumption-concentration-independence}.
We have
\begin{multline*}
\p( F(Z_{L+1}) - F(Z_L)| \ge 2t)  \le \p(|F(Z_{L+1}) - \E_L F(Z_{L+1})| \ge t) \\
+ \p(|\E_L F(Z_{L+1}) - \E F(Z_{L+1})|\ge t).
\end{multline*}

Note that for fixed values of $X,(W_\ell,B_\ell))_{\ell < L}$, the function $f \colon \R^{n_{L+1}n_L+n_{L+1}}\to \R$, given by
\begin{displaymath}
f(w_L,b_L) = F\Big(\frac{1}{\sqrt{n_{L+1}}}\sigma_L(w_L Z_L + b_L\1_n^T)\Big),
\end{displaymath}
has the Lipschitz constant bounded from above by
\begin{displaymath}
\frac{1}{\sqrt{n_{L+1}}}\lambda_\sigma(\|Z_{L}\|_{op} + \sqrt{n}).
\end{displaymath}

By the concentration property of $(W_L,B_L)$ and Lemma \ref{le:concentration-around-mean}, we therefore have $\p_{\le L-1}$ almost surely,
\begin{displaymath}
  \p_L(|F(Z_{L+1}) - \E_L F(Z_{L+1})| \ge t ) \le 2\exp\Big(-\frac{t^2}{8\rho^2 \bigl(\frac{1}{\sqrt{n_{L+1}}}\lambda_\sigma(\|Z_{L}\|_{op} + \sqrt{n})\bigr)^2}\Big),
\end{displaymath}
which, via integration by parts, and Lemma \ref{le:moments} implies that for $p \ge 1$,
\begin{align}\label{eq:conditional-moment-bound}
  \E_L |F(Z_{L+1}) - \E_L F(Z_{L+1})|^p \le D^p p^{p/2}\rho^p \frac{1}{n_{L+1}^{p/2}}\lambda_\sigma^p(\|Z_{L}\|_{op} + \sqrt{n})^p.
\end{align}
By the induction assumption $\E \|Z_{L}\|_{op} \le \E\|Z_L\|_{HS} \le C_L{\sqrt{n}}$ and
\begin{displaymath}
  \p_{\le L-1}(\|Z_L\|_{op} \ge \E \|Z_L\|_{op} + t) \le 2\exp\Big(-c_L\min_{1\le k \le L+1} (n^{(k-1)/2}t)^{2/k}\Big).
\end{displaymath}
Again, integrating by parts and using Lemma \ref{le:moments}, we obtain that for $p \ge 1$,
\begin{align*}
  \E_{\le L-1} \|Z_L\|_{op}^p \le&  2^p (\E \|Z_L\|_{op})^p + \E(\|Z_L\|_{op}-\E\|Z_L\|_{op})_+^p) \\
  \le& (2C_L)^p n^{p/2} + 2^{p}\int_0^\infty pt^{p-1}\p_{\le L-1}(\|Z_L\|_{op} \ge \E \|Z_L\|_{op} + t)dt\\
  \le& (2C_L)^p n^{p/2} + D'(L,\rho,\lambda_\sigma,\tau,\alpha,\eta)^p  \sum_{k=1}^{L+1} \frac{p^{kp/2}}{n^{(k-1)p/2}}.
\end{align*}
Combining this with \eqref{eq:conditional-moment-bound} and using the Fubini theorem, we obtain that for $p \ge 1$,
\begin{align*}
\E & |F(Z_{L+1}) - \E_L F(Z_{L+1})|^p  \le 2^pD^p \frac{n^{p/2}}{n_{L+1}^{p/2}}\rho^p p^{p/2} + 2^pD^p \rho^p n_{L+1}^{-p/2}p^{p/2} \E_{\le L-1} \|Z\|_{op}^p\\
\le& (2D \sqrt{n/n_{L+1}}\rho)^{p} p^{p/2} + (4D  C_L \rho \sqrt{n/n_{L+1}})^p p^{p/2} \\
&+(2DD'(L,\rho,\lambda_\sigma,\tau,\alpha,\eta)\rho)^p \sum_{k=1}^{L+1} \frac{p^{(k+1)p/2}}{n_{L+1}^{p/2}n^{(k-1)p/2}}\\
&\le  D''(L,\rho,\lambda_\sigma,\tau,\alpha,\eta)^p \sum_{k=1}^{L+2} \frac{p^{kp/2}}{n^{(k-1)p/2}},
\end{align*}
where in the last inequality we used the estimate $n/n_{L+1} \le \alpha$. Moreover, an inspection of the above estimates, taking into account the monotonicity of the optimal constants $c_L^{opt},C_L^{opt}$ with respect to appropriate parameters, reveals that for fixed $L,\rho,\lambda_\sigma,\tau, \eta$, we can choose the constant $D''(L,\rho,\lambda_\sigma,\tau,\alpha,\eta)$ in such a way, that
\begin{displaymath}
  \lim_{\alpha \to 0} D''(L,\rho,\lambda_\sigma,\tau,\alpha,\eta) = 0.
\end{displaymath}

Taking the $p$-th root and using the subadditivity of $s \mapsto s^{1/p}$ for $p \ge 1$, we get
\begin{align*}
  \|F(Z_{L+1}) - \E_L F(Z_{L+1})\|_p \le D''(L,\rho,\lambda_\sigma,\tau,\alpha,\eta) \sum_{k=1}^{L+2} \frac{p^{k/2}}{n^{(k-1)/2}},
\end{align*}
which by the Chebyshev's inequality implies that for all $p \ge 1$,
\begin{multline*}
  \p\Big(|F(Z_{L+1}) - \E_L F(Z_{L+1})| \ge D''(L,\rho,\lambda_\sigma,\tau,\alpha,\eta) e \sum_{k=1}^{L+2} \frac{p^{k/2}}{n^{(k-1)/2}}\Big)\\
  \le \p(|F(Z_{L+1}) - \E_L F(Z_{L+1})| \ge e \|F(Z_{L+1}) - \E_L F(Z_{L+1})\|_p) \le e^{-p}
\end{multline*}
Since for $p < 1$ , we have $e^{1-p} \ge 1$, we see that for all $p>0$, the left-hand side of the above inequality does not exceed $e^{1-p}$.

By a change of variables $p = \min_{1\le k \le L+2} (\frac{n^{(k-1)/2}t}{e(L+2)D''(L,\rho,\lambda_\sigma,\tau,\alpha,\eta)})^{2/k} $ and an adjustment of constants, this implies that
\begin{align}\label{eq:first-summand-tail}
\p(|F(Z_{L+1}) - \E_L F(Z_{L+1})| \ge t) \le 2\exp\Big(-\widetilde{c}_{L+1}\min_{1\le k \le L+2} (n^{(k-1)/2}t)^{2/k}\Big).
\end{align}
where for fixed $L,\rho,\lambda_\sigma,\tau, \eta$,
\begin{align}\label{eq:indeed-it-does-independence}
  \lim_{\alpha \to 0} \widetilde{c}_{L+1} = \infty
\end{align}
(note that if we choose $\widetilde{c}_{L+1}$ to be the best constant, such that \eqref{eq:first-summand-tail} holds for all $\XX$ satisfying the Assumption A3, then by its monotonicity with respect to $\eta$, the convergence is in fact uniform in $\eta \ge 1$, when $L,\rho,\lambda_\sigma,\tau$ are fixed).

Define now $g \colon \R^{n_{L}\times n}\to \R$ as
\begin{align}\label{eq:def-of-g}
  g(z) = \E_{L} F \Big( \frac{1}{\sqrt{n_{L+1}}}\sigma_L(W_L z + B_L\1_n^T)\Big).
\end{align}
For $z,z' \in \R^{n_{L}\times n}$ we have
\begin{align*}
  |g(z) - g(z')| &\le \E_{L}  \Big|F \bigl( \frac{1}{\sqrt{n_{L+1}}}\sigma_L(W_L z + B_L\1_n^T)\bigr) - F \bigl( \frac{1}{\sqrt{n_{L+1}}}\sigma_L(W_L z' + B_L\1_n^T)\bigr)\Big| \\
  &\le \E_L \frac{1}{\sqrt{n_{L+1}}} \lambda_\sigma \|W_{L}(z-z')\|_{HS} \\
  &\le \frac{1}{\sqrt{n_{L+1}}} \lambda_\sigma (\E_L \|W_{L}(z-z')\|_{HS}^2)^{1/2}\\
  & = \frac{1}{\sqrt{n_{L+1}}} \lambda_\sigma \Big(\sum_{i=1}^{n_{L+1}}\sum_{j=1}^n \E_L \Big(\sum_{k=1}^{n_L} W_L(i,k)(z_{k,j}-z'_{k,j})\Big)^2\Big)^{1/2}.
\end{align*}

Observe that $\E_{L} \sum_{k=1}^{n_L} W_L(i,k)(z_{k,j}-z'_{k,j}) = 0$, moreover the function
\begin{displaymath}
\R^{n_{L+1}\times n_L} \ni w \mapsto \sum_{k=1}^{n_L} w_L(i,k)(z_{k,j}-z'_{k,j})
\end{displaymath}
is $\|(z_{k,j}-z'_{k,j})_{k\le n_L}\|_2$--Lipschitz with respect to the Hilbert--Schmidt norm. Thus, by the concentration assumptions on $W_L$, Lemma \ref{le:concentration-around-mean} and integration by parts,
\begin{align}\label{eq:W-integrated}
  \E \Big(\sum_{k=1}^{n_L} W_L(i,k)(z_{k,j}-z'_{k,j})\Big)^2 \le \widetilde{D}^2\rho^2 \sum_{k=1}^{n_L}(z_{k,j}-z'_{k,j})^2
\end{align}
and as a consequence
\begin{align}\label{eq:g-is-Lipschitz}
  |g(z) - g(z')| \le \frac{1}{\sqrt{n_{L+1}}} \widetilde{D} \rho \lambda_\sigma \Big(\sum_{i=1}^{n_{L+1}}\sum_{j=1}^n \sum_{k=1}^{n_\ell}(z_{k,j}-z'_{k,j})^2\Big)^{1/2} = \widetilde{D}\rho\lambda_\sigma \|z-z'\|_{HS}.
\end{align}

Thus, $g$ is $\widetilde{D}\rho\lambda_\sigma$--Lipschitz, and by the induction assumption we get
\begin{multline*}
  \p(|\E_L F(Z_{L+1}) - \E F(Z_{L+1})|\ge t) = \p(|g(Z_L) - \E g(Z_L)|\ge t) \\
  \le 2\exp\Big(- c_L\min_{1\le k \le L+1} (n^{(k-1)/2}t)^{2/k}\Big),
\end{multline*}
where for fixed $L,\rho,\lambda_\sigma,\tau$, $\lim_{\alpha\to 0,\eta\to \infty} c_L = \infty$.

In combination with \eqref{eq:first-summand-tail} and \eqref{eq:indeed-it-does-independence} this proves \eqref{eq:induction-assumption-concentration-independence} and \eqref{eq:it-goes-tozero-induction-independence} with $L$ replaced by $L+1$.

\medskip

It remains to show that $\|Z_{L+1}\|_{HS}\le C_{L+1}\sqrt{n}$ for some constant $C_{L+1}$. By the Lipschitz property of $\sigma_L$ and the bound $|\sigma_L(0)| \le \tau$, we have
\begin{align*}
\E \|Z_{L+1}\|_{HS} & \le \sqrt{n}\tau + \frac{1}{\sqrt{n_{L+1}}} \E \|W_{L} Z_{\ell} + B_L\1_n^T\|_{HS} \\
&\le \sqrt{n}\tau + \frac{1}{\sqrt{n_{L+1}}} \E_{\le L-1}(\E_L \|W_{L} Z_{L}\|_{HS}^2)^{1/2} + \frac{\sqrt{n}}{\sqrt{n_{L+1}}}(\E\|B_L\|_2^2)^{1/2}.
\end{align*}
Using $\E B_{L} = 0$, concentration and integration by parts, we see that the last term on the right-hand side does not exceed $\widehat{D}\rho\sqrt{n}$.
As for the second term, the inequality \eqref{eq:W-integrated} applied with $z = Z_{L}$ and $z' = 0$ gives
\begin{displaymath}
  \E_L \|W_{L} Z_{L}\|_{HS}^2 \le \widetilde{D}^2\rho^2n_{L+1}\|Z_{L}\|_{HS}^2,
\end{displaymath}
so
\begin{displaymath}
  \frac{1}{\sqrt{n_{L+1}}} \E_{\le L-1}(\E_L \|W_{L} Z_{L}\|_{HS}^2) \le \widetilde{D}\rho \E\|Z_{L}\|_{HS} \le \widetilde{D}\rho C_{L}\sqrt{n},
\end{displaymath}
where in the last inequality we used the induction assumption.

This ends the proof of \eqref{eq:induction-assumption-HS-independence} and thus of the first inequality of the theorem.
The second inequality follows from the first one in the same way as in the proof of Theorem \ref{thm:unbounded-activation}.

\end{proof}
\subsection{Additional symmetry}

\begin{proof}[Proof of Theorem \ref{thm:unbounded-activation-symmetry-independence} and Proposition \ref{prop:unbounded-activation-symmetry-independence}]

The argument will be similar to the proof of Theorem \ref{thm:unbounded-activation-independence}, but thanks to the additional symmetry and the lack of bias $B_\ell$ imposed by Assumption A4, the matrices $Z_\ell$ are now centered, which will allow us to directly bound their expected operator norms by $\mathcal{O}(1+\sqrt{n/n_0})$, instead of comparing it with the Hilbert--Schmidt norm, which is typically of the order of $\sqrt{n}$. As a result we will obtain dependence on $\min(n,n_0)$ in the tail bound which improves with the growing sample size and dimension.

We will prove by induction on $L$, that there exist positive constants $c_L, C_L$, depending only on $L, \rho,\lambda_\sigma,\alpha$, such that for every random vector $\XX$, satisfying the Assumption A4,
\begin{align}\label{eq:operator-norm-induction}
\|Z_L\|_{op} \le C_L\Big(1+ \sqrt{\frac{n}{n_0}}\Big)
\end{align}
and for every $F\colon \R^{n_L\times n} \to \R$, $1$-Lipschitz with respect to the Hilbert--Schmidt norm, and any $t > 0$,

\begin{align}\label{eq:concentration-improved-induction}
\p(|F(Z_L) - \E F(Z_L)| \ge t) \le 2\exp\Big(-c_L\min(n,n_0)\min(t^2,t^{2/(L+1)})\Big).
\end{align}

As in the proof of Theorem \ref{thm:unbounded-activation-independence}, without loss of generality, we will assume that $X,W_\ell$, $\ell=0,\ldots,L$ are defined as coordinates on a product probability space. We will denote by $\E_\ell$ integration with respect to $W_\ell$ (with the other coordinates fixed) and by $\E_{\le \ell}$ integration with respect to $X,W_i$, $i \le \ell$. By the Fubini theorem we have $\E_\ell \E_{\le \ell - 1} = \E_{\le \ell}$. We will use a similar convention concerning conditional probabilities.

In what follows, in addition to the constants $C_L$ and $c_L$, we will use constants $D(a_1,\ldots,a_m)$, $D'(a_1,\ldots,a_m)$, etc., which are allowed to depend only on the parameters $a_1,\ldots,a_m$ (in particular if the list is empty, the corresponding constant is universal).

\paragraph{{\bf 1. The case $\mathbf{L=0}$}} The inequality \eqref{eq:concentration-improved-induction} with $c_0 = 1/(8\rho^2)$ follows from the concentration assumption on $X$ and Lemma \ref{le:concentration-around-mean} (note that in fact we get an estimate with $\min(n,n_0)$ replaced by $n_0$, see Remark \ref{re:no-n}). The operator norm estimate \eqref{eq:operator-norm-induction} with $C_0 = 8\rho$ follows from Lemma \ref{le:operator-norm}.

\paragraph{{\bf 2. The induction step}} Assume that \eqref{eq:operator-norm-induction} and \eqref{eq:concentration-improved-induction} hold for some constants $c_L, C_L$ and
let $\XX = (X,(W_\ell,B_\ell)_{\ell=0,\ldots,L})$ be any vector satisfying the Assumption A3 (with $L$ replaced by $L+1$, note that $Z_{L+1}$ is a function of $\XX$).

We will start with the proof of \eqref{eq:concentration-improved-induction}. Similarly as in the proof of Theorem \ref{thm:unbounded-activation-independence}, let us decompose

\begin{multline}\label{eq:decomposition}
\p( F(Z_{L+1}) - F(Z_L)| \ge 2t)  \le \p(|F(Z_{L+1}) - \E_L F(Z_{L+1})| \ge t) \\
+ \p(|\E_L F(Z_{L+1}) - \E F(Z_{L+1})|\ge t).
\end{multline}
and consider the function $f\colon \R^{n_{L+1}\times n_L} \to \R$, given by
\begin{displaymath}
f(w_L) = F\Big(\frac{1}{\sqrt{n_{L+1}}}\sigma_\ell(w_L Z_L)\Big),
\end{displaymath}
which has the Lipschitz constant bounded from above by
\begin{displaymath}
\frac{1}{\sqrt{n_{L+1}}}\lambda_\sigma \|Z_{L}\|_{op}.
\end{displaymath}

An application of the concentration assumption on $W_L$, together with Lemma \ref{le:concentration-around-mean} gives
\begin{align*}
  \p_L(|F(Z_{L+1}) - \E_L F(Z_{L+1})| \ge t) \le 2\exp\Big(-\frac{n_{L+1} t^2}{8\rho^2\lambda_\sigma^2 \|Z_L\|_{op}^2}\Big).
\end{align*}

Moreover, by the induction assumption, for $t \ge 0$,
\begin{multline*}
p(t) :=   \p_{\le L-1}(\|Z_L\|_{op} \ge C_L(1 + \sqrt{n/n_0})+ t^{(L+1)/(L+2)}) \\
\le 2\exp\Big(-c_L \min(n,n_0) \min\Big(t^{2(L+1)/(L+2)}, t^{2/(L+2)}\Big)\Big).
\end{multline*}
Combining these two estimates with the Fubini theorem, we get

\begin{multline}\label{eq:conditionally-centered-estimate}
\p(|F(Z_{L+1}) - \E_L F(Z_{L+1})| \ge t) \\
\le 2\exp\Big(-\frac{n_{L+1} t^2}{8\rho^2\lambda_\sigma^2 C_L^2 (1 + \sqrt{n/n_0} + t^{(L+1)/(L+2)})^2}\Big) + p(t)\\
\le 2\exp\Big(\frac{n_{L+1}}{72\rho^2 \lambda_\sigma^2 C_L^2 } \min\Big(t^2,\frac{n_0}{n}t^2,t^{2/(L+2)}\Big)\Big) + p(t)\\
\le 2\exp\Big(\frac{1}{72\rho^2 \alpha \lambda_\sigma^2 C_L^2 } \min\Big( n t^2,n_0 t^2,n t^{2/(L+2)}\Big)\Big) + p(t)\\
\le 4\exp\Big(-\widetilde{c}_{L+1} \min(n,n_0)\min\Big(t^2,t^{2/(L+2)}\Big)\Big)
\end{multline}
for some constant $\widetilde{c}_{L+1}$, depending only on $L,\rho,\lambda_\sigma,\alpha$.

Passing to the second term on the right-hand side of \eqref{eq:decomposition}, define $g\colon \R^{n_{L}\times n}\to \R$ as
\begin{displaymath}
  g(z) = \E_{L} F \Big( \frac{1}{\sqrt{n_{L+1}}}\sigma_L(W_L z)\Big).
\end{displaymath}
This is exactly the same definition as in \eqref{eq:def-of-g} in the proof of Theorem \ref{thm:unbounded-activation-independence}, taking into account that now, by assumption, $B_L = 0$. It is easy to see that the estimates presented below equation \eqref{eq:def-of-g} and leading to \eqref{eq:g-is-Lipschitz} are still valid in our case, leading to
the conclusion that $g$ is $D\rho\lambda_\sigma$-Lipschitz. The induction assumption \eqref{eq:concentration-improved-induction} gives therefore
\begin{multline*}
  \p(|\E_L F(Z_{L+1}) - \E F(Z_{L+1})|\ge t) = \p(|g(Z_L) - \E g(Z_L)|\ge t) \\
  \le 2\exp\Big(- c_L \min(n,n_0)\min(t^2,t^{2/(L+1)}\Big),
\end{multline*}
which in combination with \eqref{eq:decomposition} and \eqref{eq:conditionally-centered-estimate} proves \eqref{eq:concentration-improved-induction} with $L$ replaced by $L+1$ and the constant $c_L$ replaced by some constant $c_{L+1}$, depending only on $L,\rho,\lambda_\sigma,\alpha$. Note that for $t \le 1$, $\min(t^2,t^{2/(L+1)},t^{2/(L+2)}) = t^2$, while for $t > 1$ this minimum equals $t^{2/(L+2)}$.

Let us now pass to the proof of \eqref{eq:operator-norm-induction}. We will estimate $\E_L \|Z_{L+1}\|_{op}$ in terms of $\|Z_L\|_{op}$, which will allow us to use the induction assumption.

Observe that independence, together with the symmetry of $\sigma_{L+1}$ and of the law of $W_L$, implies that
\begin{displaymath}
  \E_L Z_{L+1} = \frac{1}{\sqrt{n_{L+1}}} \E_L \sigma_L(W_LZ_L) = - \frac{1}{\sqrt{n_{L+1}}} \E_L \sigma_L(W_LZ_L) = 0.
\end{displaymath}
Moreover, conditionally on $Z_L$, the random matrix $Z_{L+1}$ is a $\frac{1}{\sqrt{n_{L+1}}}\lambda_\sigma\|Z_L\|_{op}$-Lipschitz image of $W_L$ and thus satisfies the subgaussian concentration property with constant $\frac{\rho}{\sqrt{n_L}}\lambda_\sigma\|Z_L\|_{op}$. A conditional application of Lemma \ref{le:operator-norm}, combined with the Fubini theorem, gives
\begin{displaymath}
  \E \|Z_{L+1}\|_{op} =  \E_{\le L-1} \E_L \|Z_{L+1}\|_{op} \le \frac{8\lambda_\sigma\rho}{\sqrt{n_{L+1}}} (\sqrt{n_{L+1}} + \sqrt{n}) \E \|Z_L\|_{op}.
\end{displaymath}
Using the estimate $n/n_{L+1} \le \alpha$, together with the induction assumption, we can bound the right-hand side by $8C_{L}\rho\lambda_\sigma(1+\sqrt{\alpha})(1 + \sqrt{n/n_0})$, which proves \eqref{eq:operator-norm-induction} with $L+1$ instead of $L$ and $C_{L+1} = 8C_{L}\rho\lambda_\sigma(1+\sqrt{\alpha})$. This ends the induction step and the proof of the theorem.

\end{proof}

\subsection{Non-commutative polynomials in random matrices}\label{sec:polynomials-proofs}

\begin{proof}[Proof of Theorem \ref{thm:polynomials}]

The mean zero and concentration assumptions together with Lemma \ref{le:operator-norm} imply that  for any $i\le L$, $\E \|X_i\|_{op} \le  32\sqrt{n}\rho$ (the constant 32 comes from the triangle inequality, since formally we need to look at the real and imaginary part of the matrix separately). Moreover, by the concentration assumption and Lemma \ref{le:concentration-around-mean}, $\p(\|X_i\|_{op} \ge \E\|X_i\|_{op} + t) \le 2\exp(-t^2/(8\rho^2))$ and so, by the union bound,
\begin{displaymath}
  \p(\textrm{for all $i \le L$, $\|X_i\|_{op} \le 32\sqrt{n}\rho + \sqrt{8\ln (8L)} \rho$})\ge \frac{3}{4}.
\end{displaymath}
Let $C_1$ be any constant greater than $(32+\sqrt{8\ln(8 L)})\rho$ and define

\begin{displaymath}
\mathcal{A} = \{\xx \in (\C^{n\times n})^L \colon F(P(n^{-1/2}\xx)) \le \Med F(\yY) \textrm{ and } \forall_{1\le i \le L} \|X_i\|_{op} \le C_1\sqrt{n}\}.
\end{displaymath}
We have $32\sqrt{n}\rho + \sqrt{\ln (8L)} \rho \le C_1 \sqrt{n}$ and so the above tail bound and the definition of the median imply that
\begin{displaymath}
  \p(\xX \in \mathcal{A}) \ge 1/4.
\end{displaymath}
Thus, by Lemma \ref{le:enlargement} for any $t \ge 0$, $\p(\xX \notin A_t) \le 8\exp(-t^2/4\rho^2)$.

For any $\xx = (x_1,\ldots,x_L) \in \mathcal{A}$ such that $\|\xx - \xX\|_2 \le t$ we have
\begin{align}\label{eq:polynomial-upper-bound}
  F(\yY) \le F(\xx) + \|P(n^{-1/2} \xX)  - P(n^{-1/2} \xx)\|_{HS} \le \Med F(\yY) + \|P(n^{-1/2} \xX)  - P(n^{-1/2} \xx)\|_{HS}.
\end{align}
Let $\yy = (y_1,\ldots,y_L) := \xX - \xx$. We have
\begin{displaymath}
  \|P(n^{-1/2} \xX)  - P(n^{-1/2} \xx)\|_{HS} \\
  \le \sum_{k=1}^d \sum_{\ii \in [L]^k} n^{-k/2}\sum_{\varepsilon \in \{1,\ast\}^k} B_{k,\ii,\varepsilon},
\end{displaymath}
where
\begin{displaymath}
  B_{k,\ii,\varepsilon} = \Big\| A_{k,\ii,\varepsilon}^{(0)}X_{i_1}^{\varepsilon_1}A_{k,\ii,\varepsilon}^{(1)} X_{i_2}^{\varepsilon_2}\cdots A_{k,\ii,\varepsilon}^{(k-1)}X_{i_k}^{\varepsilon_k} A_{k,\ii,\varepsilon}^{(k)} -  A_{k,\ii,\varepsilon}^{(0)}x_{i_1}^{\varepsilon_1}A_{k,\ii,\varepsilon}^{(1)} x_{i_2}^{\varepsilon_2}\cdots A_{k,\ii,\varepsilon}^{(k-1)}x_{i_k}^{\varepsilon_k} A_{k,\ii,\varepsilon}^{(k)}\Big\|_{HS}
\end{displaymath}
Consider fixed $k$, $\ii = (i_1,\ldots,i_k)$ and $\varepsilon = (\varepsilon_1,\ldots,\varepsilon_k)$. Writing $X_{i_k} = x_{i_k} + y_{i_k}$, where $\yy = (y_1,\ldots,y_L)$ satisfies $\|\yy\|_2 \le t$,  and expanding the product in the definition of $B_{k,\ii,\varepsilon}$, we obtain
\begin{displaymath}
  B_{k,\ii,\varepsilon} \le
  \sum_{\emptyset \neq I \subseteq [k]} \Big\| \prod_{j=1}^d (x_{i_j}\ind{j\notin I} + y_{i_j}\ind{j \in I})A_{k,\ii,\varepsilon}^{(j)}\|_{HS}.
\end{displaymath}
For any square matrices $a_1,\ldots,a_M$ and nonempty $J \subseteq [M]$, we have
\begin{displaymath}
\|a_1\cdots a_M\|_{HS} \le \prod_{j\in J} \|a_j\|_{HS}\prod_{j\notin J}\|a_j\|_{op},
\end{displaymath}
 therefore
using the assumption $\|A_{k,\ii,\varepsilon}^{(j)}\|_{op} \le 1$, we obtain that for any nonempty $I \subseteq [k]$,
\begin{multline*}
\Big\| \prod_{j=1}^d (x_{i_j}\ind{j\notin I} + y_{i_j}\ind{j \in I})A_{k,\ii,\varepsilon}^{(j)}\Big\|_{HS}  \le \prod_{j\in I} \|y_{i_j}\|_{HS} \prod_{j\notin I}\|x_{i_j}\|_{op} \le C_1^{k- |I|} n^{(k -|I|)/2} t^{|I|} \\
= C_1^{k - |I|} n^{k/2} \Big(\frac{t}{n^{1/2}}\Big)^{|I|}.
\end{multline*}
Thus for some constants $C_3,C_4$, depending only on $d,L,\rho$,
\begin{displaymath}
  \|P(n^{-1/2} \xX)  - P(n^{-1/2} \xx)\|_{HS} \le C_3 \sum_{i=1}^d  \Big(\frac{t}{n^{1/2}}\Big)^i
  \le C_4 \Big(\frac{t}{\sqrt{n}} + \Big(\frac{t}{\sqrt{n}}\Big)^{d}\Big).
\end{displaymath}
Going back to \eqref{eq:polynomial-upper-bound} we obtain that
\begin{displaymath}
\p\Big(F(\yY) \ge \Med F(\yY) + C_4\Big(\frac{t}{\sqrt{n}} + \Big(\frac{t}{\sqrt{n}}\Big)^{d}\Big)\Big) \le 8\exp(-t^2/4\rho^2),
\end{displaymath}
which implies that for some constant $c$, depending only on $d,L,\rho$ and all $t \ge 0$
\begin{displaymath}
  \p(F(\yY) \ge \Med F(\yY) + t)\le 8\exp\Big(-cn \min(t^2, t^{2/d})\Big).
\end{displaymath}
Applying this inequality to $-F$, using the union bound and adjusting the constants, we obtain the statement of the theorem.
\end{proof}

\begin{remark}\label{re:trace}
In Remark \ref{le:polynomials} we mentioned that in the case of $f(x) = x$, i.e. when dealing with the trace of the matrix, one can replace $n^{1-1/d}$ in \eqref{eq:les-polynomials} with $n$. This can be done, by observing that thanks to the linearity of the trace, one does not need to pass immediately to the Hilbert-Schmidt norm, but can analyse the quantities of the form
$\tr \prod_{j=1}^d (x_{i_j}\ind{j\notin I} + y_{i_j}\ind{j \in I})A_{k,\ii,\varepsilon}^{(j)}$. For $|I| \ge 2$, using the cyclic invariance of the trace together with the estimate $\tr AB \le \|A\|_{HS}\|B\|_{HS}$, we can thus estimate this quantity by $\prod_{j\in I} \|y_{i_j}\|_{HS} \prod_{j\notin I}\|x_{i_j}\|_{op}$, which gives an improvement by a factor $\sqrt{n}$ with respect to the quantity appearing in the above proof (note that the trace is $\sqrt{n}$-Lipschitz with respect to the Hilbert--Schmidt norm). For $|I| = 1$ this improvement does not appear and we still need to rely on the comparison between the trace and the Hilbert--Schmidt norm to be able to combine the bounds on $\|y_i\|_{HS}$ and $\|x_i\|_{op}$. This explains why in the case of the trace the subgaussian coefficient of \eqref{eq:les-polynomials} remains independent of $n$, while the other coefficient can be improved (which agrees with the inequality by Meckes and Szarek \cite{MR2869165}).
\end{remark}
\section{Proofs of results from Section \ref{sec:applications}}\label{sec:applications-proofs}

We will now prove corollaries to our main inequalities related to convergence of empirical spectral measures and concentration for the conjugate kernels.

\begin{proof}[Proof of Corollary \ref{cor:distance-bound}]
By Propositions \ref{prop:unbounded-activation} and \ref{prop:unbounded-activation-independence}, we have $\E \|Z_L\|_{HS} \le C\sqrt{n^{(m)}}$ for some constant $C$, independent of $m$. By Theorems \ref{thm:unbounded-activation} and \ref{thm:unbounded-activation-independence}, integration by parts and the triangle inequality we obtain for some other constants $C',C''$, also independent of $m$,
\begin{displaymath}
  \E \|Z_L\|_{HS}^2 \le 2 (\E \|Z_L\|_{HS})^2 + C' \le C'' n.
\end{displaymath}
Since
\begin{displaymath}
  \frac{1}{n} \E \|Z_L\|_{HS}^2 = \int_{\R} x^2 \E \mu_{A_m}(dx),
\end{displaymath}
this implies that the sequence $(\E \mu_{A_m})_{m\ge 1}$ is tight. We can therefore assume without loss of generality that this sequence converges weakly to some probability measure $\mu$. Indeed, if one could find a counterexample to Corollary \ref{cor:distance-bound}, by passing to a subsequence one could find a counterexample satisfying this additional assumption.

It is thus enough to show that $d(\mu_{A_m},\mu) \to 0$ in probability. Let $\widetilde{\mu}$ be the image of $\mu$ under the continuous map $x \mapsto \sqrt{x}$. The image of $\mu_{A_m}$ under this map equals to $\mu_{\sqrt{A_m}}$, and the image of $\E \mu_{A_m}$ is $\E \mu_{\sqrt{A_m}}$. In particular $\E \mu_{\sqrt{A_m}}$ converges weakly to $\widetilde{\mu}$, moreover the desired convergence  $d(\mu_{A_m},\mu) \to 0$ is equivalent to the weak in probability convergence of $\mu_{\sqrt{A_m}}$ to $\widetilde{\mu}$, which can be rephrased as the convergence in probability
\begin{displaymath}
  \int_\R f d\mu_{\sqrt{A_m}} \to \E \int_\R f d\widetilde{\mu}
\end{displaymath}
for all bounded Lipschitz functions $f$ or equivalently, replacing $\widetilde{\mu}$ by $\E \mu_{\sqrt{A_m}}$ and passing to the notation introduced in \eqref{eq:les} as
\begin{displaymath}
  \frac{1}{n}(S_f(Z_L^{(m)}) - \E S_f(Z_L^{(m)})) \to 0
\end{displaymath}
in probability (see, e.g., \cite[Lemma A.1]{MR4302287}).
The last convergence is an immediate consequence of Corollary \ref{cor:les-unbounded-activation}. This ends the proof.
\end{proof}

\begin{proof}[Proof of Corollary \ref{cor:les-conjugate-kernel}]
Let $\Gamma=(\gamma_1,\ldots,\gamma_n)$, where $\gamma_i$ are the singular values of $Z_L$ arranged in non-increasing order. Then $Y = (\gamma_1^2,\ldots,\gamma_n^2)$. Moreover, by the Hoffman-Wielandt inequality (Theorem \ref{thm:HW}), $\Gamma$ is a 1-Lipschitz image of $Z_L$ and so, by Theorem \ref{thm:unbounded-activation-symmetry-independence} and Lemma \ref{le:concentration-equivalence}, for any $1$-Lipschitz $G\colon \R^N \to \R$, and $t > 0$,
\begin{displaymath}
  \p(|F(\Gamma) - \Med F(\Gamma)| \ge t) \le 2\exp(-d n\min(t^2,t)).
\end{displaymath}
where $d$ depends only on $L,\rho,\lambda_\sigma,\alpha$. By Proposition \ref{prop:unbounded-activation-symmetry-independence}, for some constant $C$, depending only on $L,\rho,\lambda_\sigma,\alpha$, $\p(\|Z_L\|_{op} \ge C) \le 1/4$. Since $\|Z_L\| = \max_{i\le n} \gamma_i$, we conclude that $\mathcal{A} = \{x = (x_i)_{i\le n}\colon F(x_1^2,\ldots,x_n^2) \le \Med F(Y), \max_{i\le n} x_i  \le C\}$ satisfies $\p(\Gamma \in \mathcal{A}) \ge 1/4$. The following lines resemble closely the proof of Lemma \ref{le:enlargement}.

Set
\begin{displaymath}
t_0 = \max\Big(\Big(\frac{\ln 8}{dn}\Big)^{1/2},\Big(\frac{\ln 8}{dn}\Big)^{(L+1)/2}\Big).
\end{displaymath}
Then for any $t > t_0$, $\p(\Gamma \in \mathcal{A}_{t}) \ge 1/2$. Indeed, otherwise, we would have $\p(\Gamma \in \mathcal{A}_{t}^c) \ge 1/2$ and so, by concentration
\begin{displaymath}
  \p(\Gamma \in (\mathcal{A}_{t}^c)_{t}) \ge 1 - 2\exp(-dn\min(t,t^{2/(L+1)}) > 3/4,
\end{displaymath}
which is a contradiction, since $\mathcal{A} \cap ((\mathcal{A}_{t})^c)_{t} = \emptyset$.

Thus, for any $s > 2t_0$,
\begin{multline*}
  \p(\Gamma \notin \mathcal{A}_{s}) \le \p(\Gamma \notin (\mathcal{A}_{s/2})_{s/2}) \le 2\exp(- d n \min ((s/2)^2,(s/2)^{2/(L+1)})) \\
  \le 2\exp(2^{-1}d n \min(s^2,s^{2/(L+1)})),
\end{multline*}
while for $s < 2t_0$ we trivially have
\begin{displaymath}
  \p(\Gamma \notin \mathcal{A}_{s}) \le 1 \le 8 \exp(-4^{-1} dn \min(s^2,s^{2/(L+1)})).
\end{displaymath}
Thus, for all $s > 0$, we have $\p(Y \notin \mathcal{A}_s) \le 8\exp(-d' n \min(s^2,s^{2/(L+1)}))$.

If $\Gamma \in \mathcal{A}_s$, then there exists $x = (x_i)_{i\le n} \in \mathcal{A}$, such that $\|x - \Gamma\| \le s$. As a consequence
\begin{multline*}
  F(Y) \le  F(x_1^2\ldots,x_n^2) + \Big(\sum_{i=1}^n |\gamma_i^2 - x_i^2|^2\Big)^{1/2} \le \Med F(Y) + \Big(\sum_{i=1}^n |\gamma_i^2 - x_i^2|^2\Big)^{1/2}\\
  \le \Med F(Y) + \Big(2\sum_{i=1}^n (\gamma_i - x_i)^4 + 8\sum_{i=1}^n x_i^2 (\gamma_i - x_i)^2\Big)^{1/2}\\
  \le \Med F(Y) + \sqrt{2}\Big(\sum_{i=1}^n (\gamma_i - x_i)^4\Big)^{1/2} + 2\sqrt{2}\|\gamma - x\|_2 \max_{i\le n} x_i \\
  \le \Med F(Y) + \sqrt{2} s^2 + 2\sqrt{2}C s,
\end{multline*}
where in the last inequality we used the definition of $\mathcal{A}$. Thus,
\begin{displaymath}
  \p(F(Y) \ge \Med F(Y) + \sqrt{2}s^2 + 2\sqrt{2}C s) \le \p(\Gamma \notin \mathcal{A}_s) \le 8\exp( -d n \min(s^2, s^{2/(L+1)})).
\end{displaymath}
Substituting $s = \min(\sqrt{t/2\sqrt{2}},t/4\sqrt{2}C)$ and adjusting the constants we obtain
\begin{multline*}
  \p(F(Y) \ge \Med F(Y) + t) \le \p(F(Y) \ge \Med F(Y) + 2\sqrt{2}s^2 + 2\sqrt{2}Cs) \\
  \le 8\exp( - \widetilde{c} n\min(t^2,t^{1/(L+1)}))
\end{multline*}
for some positive constant $\widetilde{c}$, depending only on $L,\rho, \lambda_\sigma, \alpha$. The corollary with the expectation replaced by the median follows now by applying the above inequality to $-F$, using the union bound and adjusting the constants. To pass from the median to the expectation one uses Lemma \ref{le:concentration-equivalence}.
\end{proof}

Let us conclude this section with the proof of Proposition \ref{prop:Wasserstein-convergence} and Corollaries \ref{cor:Wasserstein-conjugatekernel}, \ref{cor:Wasserstein-polynomials}.

\begin{proof}[Proof of Proposition \ref{prop:Wasserstein-convergence}]
Let  $Lip_0$ be the class of Lipschitz functions $f\colon \R\to \R$, such that $f(0) = 0$, equipped with the norm $\|\cdot\|_{Lip}$ defined as the optimal Lipschitz constant.
Denote by $B(Lip_0)$ the unit ball of $Lip_0$.
Using \eqref{eq:Wasserstein-duality} we represent
\begin{displaymath}
  \mathcal{W}_1(\mu_A,\E \mu_A) = \sup_{f\in B(Lip_0)} |X_f|,
\end{displaymath}
where $X_f = \int_\R fd\mu_A - \E \int_\R  fd\mu_A = n^{-1}(S_f(A) -\E S_f(A))$, where $S_f(A) = n\int_R f d\mu_A$.

Observe that since for $f \in B(Lip_0)$, $|f(x)| \le |x|$, the condition $\mu_A(\{x\in \R\colon  |f(x)| > t\}) > 0$ implies that $\mu_A(\R \setminus [-t,t]) > 0$, which in turn gives $\|A\|_{op} > t$ and $\mu_A(\{x\in \R\colon  |f(x)| > t\}) \le  1 = \ind{\|A\|_{op} > t}$. Thus for all $f \in Lip_0$ and $R>0$,
\begin{multline*}
\int_\R |f| \ind{\R\setminus [-R,R]} d\mu_A = \int_0^\infty \mu_A( |f| \ind{\R\setminus [-R,R]} > t) dt \le R \mu_A(\R \setminus [-R,R]) + \int_{R}^\infty \ind{\|A\|_{op} > t} dt\\
\le R \ind{\|A\|_{op} > R}  + \int_{R}^\infty \ind{\|A\|_{op} > t} dt.
\end{multline*}
As a consequence
\begin{multline*}
\E \sup_{f \in Lip_0} \Big|\int f \ind{\R\setminus [-R,R]} d\mu_A - \E \int f \ind{\R\setminus [-R,R]} d\mu_A\Big| \le 2 \E \sup_{f \in Lip_0} \int_R |f| \ind{\R\setminus [-R,R]} d\mu_A \\
\le 2R \p(\|A\|_{op} \ge R) + \int_R^\infty \p(\|A\|_{op} > t)dt.
\end{multline*}
Setting $R = 2D$, and using the assumption $\E \|A\|_{op} \le D$, together with the fact that the operator norm is a $1$-Lipschitz function of $A$ and the concentration assumption, we thus get
\begin{multline*}
\E \sup_{f \in Lip_0} \Big|\int f \ind{\R\setminus [-R,R]} d\mu_A - \E \int f \ind{\R\setminus [-R,R]} d\mu_A\Big| \\
\le 4R\exp(-an \min ((R/2)^2,(R/2)^r)) + 4\int_{R}^\infty \exp(-an \min((t/2)^2,(t/2)^r))dt \\
\le C_{a,r} \exp(- c_{a,r} n \min(D^2,D^r)).
\end{multline*}

Thus denoting by $Lip_0(2D)$ the set of $1$-Lipschitz functions such that $f(0)=0$ and $f$ is constant on $\R\setminus [-2D,2D]$, we  obtain
\begin{displaymath}
  \E \mathcal{W}_1(\mu_A,\E \mu_A) = \E \sup_{f\in B(Lip_0)} |X_f| \le \E \sup_{f\in B(Lip_0(2D))} |X_f| + 2 C_{a,r} \exp(- c_{a,r} n \min(D^2,D^r)).
\end{displaymath}

Fix now a positive integer $M$ and denote by $\mathcal{X}$ the linear space of all functions $f \in Lip_0(2D)$, which are affine on each interval of the form $(2Dk/M,2D(k+1)/M)$, and $[-2D(k+1)/M,-2Dk/M)$, $k = 0,\ldots,M-1$.  Thus $\mathcal{X}$ is a linear space of dimension $2M$. Moreover, for any $f \in B(Lip_0(2D))$, there exists $\widetilde{f} \in B(\mathcal{X})$, such that $\|f - \widetilde{f}\|_\infty \le 2D/M$, thus
\begin{align}\label{eq:intermediate-Wasserstein}
  \E \mathcal{W}_1(\mu_A,\E \mu_A) \le \E \sup_{f\in B(\mathcal{X})} |X_f| + 2 C_{a,r} \exp(- c_{a,r} n \min(D^2,D^r)) + 4D/M.
\end{align}

Denote by $N(B(\mathcal{X}),\varepsilon)$ the minimal number of closed balls (in the metric $\|\cdot\|_{Lip}$) needed to cover $B(\mathcal{X})$. From the standard volumetric estimate (see, e.g., \cite[Lemma 2.6]{MR856576}) it follows that for $\varepsilon < 1$, $N(B(\mathcal{X}),\varepsilon)) \le (3/\varepsilon)^{2M}$.

By the concentration assumption, for $f, g \in Lip_0$ and $t \ge 0$, we have
\begin{multline*}
  \p(|X_f - X_g| \ge t) = \p( |S_{f-g}(A) - \E S_{f-g}(A)| \ge n t) \\
  \le 2\exp\Big(-a \min\Big(n^2 \frac{t^2}{\|f-g\|_{Lip}^2},n^{1+r/2} \frac{t^{r}}{\|f-g\|_{Lip}^r}\Big)\Big),
\end{multline*}
where we used that $S_f(A)$ is a $\sqrt{n}$--Lipschitz function of $A$ (see Theorem \ref{thm:HW}).  Thus, by a generalization of Dudley's entropy bound, we get
\begin{multline*}
  \E \sup_{f \in B(\mathcal{X})} |X_f| \le C'_r \Big( \frac{1}{\sqrt{a}n} \int_0^1 \sqrt{\ln N(B(\mathcal{X}),\varepsilon)}d\varepsilon + \frac{1}{a^{1/r}n^{1/r+1/2}}\int_0^1 (\ln N(B(\mathcal{X}),\varepsilon))^{1/r}d\varepsilon\Big)\\
  \le C_{r}'' \Big(\frac{\sqrt{M}}{\sqrt{a}n}+ \frac{M^{1/r}}{a^{1/r}n^{1/r+1/2}}\Big).
\end{multline*}
The above version of Dudley's theorem is a part of folklore and we have not been able to find an exact reference. We refer to \cite{MR4381414,MR1666908,MR1385671} for general description of chaining methods, including Dudley-type bounds.

Going back to \eqref{eq:intermediate-Wasserstein}, we obtain
\begin{displaymath}
  \E \mathcal{W}_1(\mu_A,\E \mu_A)
  \le C_{a,r}'\Big(\frac{\sqrt{M}}{\sqrt{a}n}+ \frac{M^{1/r}}{a^{1/r}n^{1/r+1/2}} + \exp(- c_{a,r} n \min(D^2,D^r)) + 4D/M\Big)
  \le \frac{C_{a,r,D}}{n^{2/3}},
\end{displaymath}
where in the last line we substituted $M \simeq n^{2/3}$ and used the assumption $r \le 2$.

To finish the proof it is enough to note that by \eqref{eq:Wasserstein-duality} and Theorem \ref{thm:HW}, for any probability measure $\nu$ on $\R$, the function $A \mapsto \mathcal{W}_1(A,\nu)$ is $n^{-1/2}$-Lipschitz with respect to the Hilbert--Schmidt norm.
\end{proof}

\begin{proof}[Proof of Corollary \ref{cor:Wasserstein-conjugatekernel}] Let $\hat{A}$ be the diagonalization of $A$. Then, by Corollary \ref{cor:les-conjugate-kernel}, it satisfies the concentration assumption \eqref{eq:Wasserstein-concentration-assumption} of Proposition \ref{prop:Wasserstein-convergence} with $r = 1/(L+1)$. Since the spectra of $\hat{A}$ and $A$ are the same, it is thus enough to show that $\E\|A\|_{op} \le D$ for some $D$, depending only on  $L,\rho,\lambda_\sigma,\alpha$. Note that $\|A\|_{op} = \|Z_L\|_{op}^2$. Thus,
\begin{multline*}
  \E\|A\|_{op} = \E\|Z_L\|_{op}^2 \le 2(\E \|Z_L\|_{op})^2 + 2\E (\|Z_L\|_{op} - \E \|Z_L\|_{op})^2\\
  \le C (1+ \sqrt{n/n_0})^2 + 4\int_0^\infty t \p(\|Z_L\|_{op} \ge \E \|Z_L\|_{op} + t)dt\\
  \le C (1+ \sqrt{n/n_0})^2 + 8\int_0^\infty t \exp(-cn\min(t^2,t^{2/(L+1)})dt \le D,
\end{multline*}
where in the second inequality we used Proposition \ref{prop:unbounded-activation-symmetry-independence}, in the third inequality Theorem \ref{thm:unbounded-activation-symmetry-independence} and in the last inequality the assumption $n/n_0 \le \alpha$.
The statement follows now from Proposition \ref{prop:Wasserstein-convergence}.
\end{proof}

\begin{proof}[Proof of Corollary \ref{cor:Wasserstein-polynomials}]
The concentration assumption \eqref{eq:Wasserstein-concentration-assumption} of Proposition \ref{prop:Wasserstein-convergence} with $r = 2/d$ follows by Theorem \ref{thm:polynomials}. By Lemma \ref{le:operator-norm}, we have $\E \|X_i\|_{op} \le 32\sqrt{n}\rho$ and integration by parts, similar as in the proof of Corollary \ref{cor:Wasserstein-conjugatekernel}, shows that $\E \|n^{-1/2}X_i\|_{op}^p \le C(p,\rho)$ for any $p \ge 1$. Since the operator norm is submultiplicative, the H\"older inequality together with \eqref{eq:polynomial} and the bound $\|A_{k,\ii,\varepsilon}\|_{op} \le 1$ implies that $\E \|P(n^{-1/2}\xX)\|_{op} \le D$ for some $D$, depending only on $L,\rho,d$, which allows to end the proof by an application of Proposition \ref{prop:Wasserstein-convergence}.
\end{proof}

\appendix
\section{Additional proofs}\label{sec:appendix}
In this section we gather proofs of elementary technical facts which were used throughout the paper. Many of these facts are parts of folklore and we provide proofs just for the sake of completeness (when we were not able to locate them in the literature) or to get explicit dependence of constants on parameters.

\begin{proof}[Proof of Proposition \ref{prop:independent-concentration}]
The proof is just an adaptation of the classical martingale argument (see, e.g., \cite{MR1849347}) to the subgaussian setting. Assume without loss of generality that the random variables $X_1,\ldots,X_L$ are defined as coordinates on a product probability space and denote by $\E_{\ge i}$ integration with respect to the coordinates $X_{i},\ldots,X_{L}$ with the other coordinates fixed, while by $\E_i$, integration just with respect to $X_i$ (with the same convention for probabilities). In particular $\E_{\ge L + 1} F(X_1,\ldots,X_L) = F(X_1,\ldots,X_L)$, $\E_{\ge 1} F(X_1,\ldots,X_L) = \E F(X_1,\ldots,X_L)$. We have
\begin{displaymath}
  F(\XX) - \E F(\XX) = \sum_{i=1}^{L} \E_{\ge i+1} F(\XX) - \E_{\ge i} F(\XX).
\end{displaymath}
Observe that for any real random variable $Y$, such that for all $t > 0$, $\p(|Y| \ge t) \le 2\exp(-t^2/a^2)$, we have for all $k \in \N$,
\begin{displaymath}
  \E Y^{2k} = \int_0^\infty 2k t^{2k-1} \p(|Y|\ge t)dt \le 2 \int_0^\infty 2k t^{2k-1} e^{-t^2/a^2}dt = 2k a^k \Gamma(k) \le 2 a^{2k} k!,
\end{displaymath}
where we used Lemma \ref{le:moments}.
When $Y$ is centered and $Y'$ is an independent copy of $Y$, then for any $s \in \R$,
\begin{multline*}
  \E e^{sY} = \E e^{s(Y - \E Y')} \le \E e^{s (Y-Y')} = \sum_{k=0}^\infty \frac{s^k}{k!}\E (Y-Y')^k \\
 \le  \sum_{k=0}^\infty \frac{s^{2k}}{(2k)!}\E (Y-Y')^{2k} \le \sum_{k=0}^\infty \frac{2^{2k}s^{2k}}{(2k)!}\E Y^{2k} \le 1 + 2\sum_{k=1}^\infty \frac{2^{2k}s^{2k}}{(2k)!} a^{2k} k!\\
  \le  \sum_{k=0}^\infty \frac{2^{2k}s^{2k}}{k!} a^{2k} \le \exp(4 s^2 a^2).
  \end{multline*}

Observe that if $F$ is 1-Lipschitz on $\R^{n_1+\ldots+n_L}$, then for fixed $x_1,\ldots,x_{i-1}$ the function $x \mapsto \E_{\ge i+1} F(x_1,\ldots,x_{i-1},x,X_{i+1},\ldots,X_L)$ is 1-Lipschitz on $\R^{n_i}$ and thus by the concentration assumption on $X_i$ and Lemma \ref{le:concentration-around-mean},
\begin{displaymath}
  \p_i(\E_{\ge i+1} F(\XX) - \E_{\ge i} F(\XX)) \le 2\exp(-t^2/8\rho^2) \textrm{ a.s.}
\end{displaymath}
As a consequence, by the above calculation for $Y = \E_{\ge i+1} F(\XX)$,
\begin{displaymath}
  \E_{i} \exp(s \E_{\ge i+1} (F(\XX) - \E_{\ge i} F(\XX))) \le \exp(32 s^2 \rho^2) \textrm{ a.s.}
\end{displaymath}
Using the fact that $\E_{\ge i+1} F(\XX) - \E_{\ge i} F(\XX)$, $i \le j$, is independent of $X_{j+1},\ldots,X_L$, we can thus successively estimate,
\begin{multline*}
  \E \exp(s(F(\XX) - \E F(\XX))) = \E \prod_{i=1}^L \exp(s(\E_{\ge i+1} F(\XX) - \E_{\ge i} F(\XX))) \\
  \le \E \prod_{i=1}^{L-1} \exp(s(\E_{\ge i+1} F(\XX) - \E_{\ge i} F(\XX)))
  \E_L \exp(s(\E_{\ge L+1} F(\XX) - \E_{\ge L} F(\XX))) \\
  \le \E \prod_{i=1}^{L-1} \exp(s(\E_{\ge i+1} F(\XX) - \E_{\ge i} F(\XX))) e^{32s^2\rho^2} \le \cdots \le \exp(32L\rho^2 s^2).
\end{multline*}

Using Chebyshev's inequality we thus get
\begin{displaymath}
  \p(F(\XX) - \E F(\XX) \ge t) \le \inf_{s \ge 0} e^{32L\rho^2 s^2 - st} = \exp(- t^2/(128 L \rho^2)).
\end{displaymath}
Applying the same inequality to $-F$ and using Lemma \ref{le:concentration-equivalence} ends the proof.
\end{proof}

\begin{proof}[Proof of Lemma \ref{le:enlargement}]

The assumption \eqref{eq:concentration-assumption} is equivalent to the statement that for any $C \subseteq \R^N$, such that $\p(X \in C) \ge 1/2$, and any $t > 0$, $\p(X \notin C_t) \le 2\exp(-t^2/\rho^2)$.

Consider any $t_0$, such that $\p(X \in A_{t_0}) \le 1/2$ and set $B = (A_{t_0})^c$. Then $\p(\XX \in B) \ge 1/2$. Note that $A\cap B_{t_0} = \emptyset$. Thus $p \le \p(X \in A) \le \p(X \notin B_{t_0}) \le 2e^{-t_0^2/\rho^2}$, i.e.,
\begin{displaymath}
  t_0 \le \rho \sqrt{\ln(2/p)}.
\end{displaymath}
It follows that for any $t> t_0 > \rho\sqrt{\ln(2/p)}$, $\p(X \in  A_{t_0}) \ge 1/2$, and so
\begin{displaymath}
  \p(X \notin A_{t}) \le \p(X \notin (A_{t_0})_{t-t_0}) \le 2e^{-(t-t_0)^2/\rho^2}.
\end{displaymath}
Thus for $t > 2 \rho\sqrt{\ln(2/p)}$,
\begin{displaymath}
  \p(X \notin A_t) \le 2\exp\bigl(-(t- \rho \sqrt{\ln(2/p)})^2/\rho^2\bigr) \le 2\exp(-t^2/(4\rho^2)) \le \frac{2}{p}\exp(-t^2/(4\rho^2)),
\end{displaymath}
while for $t \le  2\rho\sqrt{\ln(2/p)}$, we have
\begin{displaymath}
  \frac{2}{p}\exp(-t^2/(4\rho^2)) \ge 1 \ge \p(X \notin A_t).
\end{displaymath}
\end{proof}

\begin{proof}[Proof of Lemma \ref{le:concentration-around-mean}]
We have
\begin{multline*}
  \E F(Y) - \Med F(Y) \le (\E F(Y) - \Med F(Y))_+ \le \E (F(Y) - \Med F(Y))_+\\
   = \int_0^\infty \p(F(Y) \ge  \Med F(Y) + t)dt
  \le \int_0^{\rho \ln 4} \frac{1}{2}dt + 2\int_{\rho \ln 4}^\infty e^{-t^2/\rho^2}dt \\
  \le \frac{\rho \ln 4}{2} + e^{- (\ln 4)^2}\cdot 2\int_{\rho \ln 4}^\infty e^{-(t-\rho \ln 4)^2/\rho^2} dt\\
  = \frac{\rho \ln 4}{2} + e^{- (\ln 4)^2} \sqrt{\pi}\rho \le \rho,
\end{multline*}
where in the second estimate we used Jensen's inequality, in the fourth one the estimate $x^2 \ge (x-a)^2 + a^2$ valid for $x \ge a \ge 0$, and in the last one a numerical calculation $\frac{\ln 4}{2}+  e^{-(\ln 4)^2}\sqrt{\pi} \simeq 0.953$. Applying the above inequality to $-F$, we get $|\E F(Y) - \Med F(Y)| \le \rho$.

Therefore, for $t \ge 2\rho$,
\begin{displaymath}
  \p(|F(Y) - \E F(Y)| \ge t) \le \p(|F(Y) - \Med F(Y)| \ge t/2) \le 2\exp(-t^2/4\rho^2).
\end{displaymath}
On the other hand for $t < 2\rho$, $2\exp(-t^2/8\rho^2) \ge 2 e^{-1/2} \ge 1$, so the tail bound of the lemma becomes trivial.
\end{proof}

\begin{proof}[Proof of Lemma \ref{le:concentration-equivalence}]
  Assume first that \eqref{eq:conc-median} holds. Then
  \begin{multline*}
    |\E F(Y) - \Med F(Y)| \le \E |F(Y) - \E F(Y)|  = \int_0^\infty \p(|F(Y) - \E F(Y)| \ge t)dt \\
    \le 2\sum_{i=1}^m \int_0^\infty e^{-a_i t^\alpha_i} dt \le C \max_{i\le m} \frac{1}{a_i^{-1/\alpha_i}}
  \end{multline*}
  for some constant $C$, depending only on $\alpha_1,\ldots,\alpha_m$.

Denote the right hand side above by $t_0$. For $t > 2t_0$, we have
\begin{displaymath}
  \p(|F(Y) - \E F(Y)| \ge t) \le \p(|F(X) - \Med F(X)| \ge t/2) \le 2 \exp(-\min_{i\le m} (2^{-\alpha_i} a_i t^{\alpha_i})).
\end{displaymath}

On the other hand, for $t < 2 t_0$, and $\widetilde{c}$ small enough, depending only on $\alpha_1,\ldots,\alpha_m$, we have
\begin{displaymath}
2\exp(-\min_{i\le m} (\widetilde{c} a_i t^{\alpha_i})) \ge 2\exp( - \max_{i\le m} (\widetilde{c} 2^{\alpha_i} C^{\alpha_i})) \ge 1 \ge \p(|F(Y) - \E F(Y)| \ge t).
\end{displaymath}
This proves that \eqref{eq:conc-median} implies \eqref{eq:conc-mean} with $b_i \ge c a_i$ for some $c$, depending only on $\alpha_1,\ldots,\alpha_m$.

Let us now assume \eqref{eq:conc-mean}. Then for $t_0 = \max_{i\le m} (b_i^{- 1/\alpha_i}  (3 \ln 2)^{1/\alpha_i})$, we have
\begin{displaymath}
\p(F(Y) \ge \E F(Y  ) + t), \p(F(Y) \le \E F(Y) - t) \le 1/4,
\end{displaymath}
which implies that $|\Med F(Y) - \E F(Y)| \le t_0$. The rest of the proof is analogous as in the first implication.
\end{proof}

\begin{proof}[Proof of Lemma \ref{le:operator-norm}]
Let $\mathcal{K}$, resp. $\mathcal{M}$ be $1/4$-nets in the unit balls of $\R^k$, resp. $\R^m$, of cardinality $9^k$, resp. $9^m$ (by standard volumetric estimates for any $\varepsilon > 0$ there exists an $\varepsilon$-net of cardinality at most $(1+2/\varepsilon)^n$, see, e.g., \cite[Lemma 2.6]{MR856576}). Then for any unit vectors $x \in \R^k$, $y \in \R^m$, if $x' \in \mathcal{K}$, $y' \in \mathcal{M}$ are such that $\|x-x'\|_2,\|y-y'\|_2 \le 1/4$,
\begin{multline*}
  |y^T U x| \le |(y - y')^T U (x-x')| + |(y')^T U (x - x')| + |(y-y')^TU x'| + |(y')^T U x'| \\
  \le \frac{9}{16}\|U\|_{op} + \max_{x'\in \mathcal{K},y'\in \mathcal{M}} |(y')^TU x'|.
\end{multline*}
Taking supremum over $x,y$, we obtain
\begin{displaymath}
  \|U\| \le \frac{16}{7} \max_{x'\in \mathcal{K},y'\in \mathcal{M}} |(y')^TU x'|.
\end{displaymath}
The first inequality of the lemma follows now by the union bound, and Lemma \ref{le:concentration-around-mean} combined with the observation that for each $y',x'$, the function $U \mapsto (y')^TU x'$ is 1-Lipschitz.

As for the bound on the expectation, for $s \ge 3$ we have
\begin{multline*}
\p(\|U\| > \frac{16}{7} s(\sqrt{k} + \sqrt{m})\rho) \le 9^{k+m}\cdot 2\exp(s^2(\sqrt{k}+\sqrt{m})^2/2)\\
\le  2\exp(2.25(k+m) - s^2(\sqrt{k}+\sqrt{m})^2/2) \le 2\exp((k+m)(4.5 - s^2)/2) \le 2\exp(-s^2/2).
\end{multline*}
and so
\begin{displaymath}
  \E\frac{7\|U\|}{16(\sqrt{k} + \sqrt{m})\rho} = \int_0^\infty \p\Big(\|U\|\ge \frac{16}{7}s(\sqrt{k} + \sqrt{m})\Big) ds
  \le 3 + \int_3^{\infty} 2\exp(-s^2/2) ds \le 3.5.
\end{displaymath}
\end{proof}
\bibliographystyle{amsplain}	
\bibliography{NN.bib}
\end{document}